\documentclass[a4paper,reqno]{amsart}
\usepackage{amsfonts}
\usepackage{amssymb}
 \usepackage{times}
 \usepackage{graphicx}

\DeclareMathOperator{\id}{id}
\newcommand{\CD}{  \mathfrak{CD}}
\DeclareMathOperator{\End}{End}
\DeclareMathOperator{\Aut}{Aut}
\DeclareMathOperator{\tr}{tr}
\DeclareMathOperator{\Mat}{Mat}
\DeclareMathOperator{\ad}{ad}
\DeclareMathOperator{\Der}{Der}
\DeclareMathOperator{\supp}{supp}
 \DeclareMathOperator{\Diag}{Diag}

\DeclareMathOperator{\Hom}{Hom}
\DeclareMathOperator{\Cent}{Cent}
\DeclareMathOperator{\PSL}{PSL}

\newcommand{\FF}{\mathbb{F}}
\newcommand{\pF}{p\mathbb{F}}
\newcommand{\K}{ \mathcal{K}}
\newcommand{\Q}{ \mathcal{Q}}
\newcommand{\C}{ \mathcal{C}}
\newcommand{\Ok}{\mathrm{Ok}}
\newcommand{\pK}{p\mathcal{K}}
\newcommand{\pQ}{p\mathcal{Q}}
\newcommand{\pC}{p\mathcal{C}}
\newcommand{\ii}{\textbf{i}}
\def\span#1{\langle #1\rangle}

\newtheorem{Definition}{Definition}
\newtheorem{Theorem}{Theorem}

\newtheorem{Remark}{Remark}
\newtheorem{Example}{Example}

\begin{document}

 \title{ Fine gradings \\on the simple Lie algebras of type $E$}

 \author{Cristina Draper${}^\star$}
  \address{Department of Applied Mathematics, University of M\'{a}laga}
    \email{cdf@uma.es}
    \thanks{${}^\star$ Supported by the Spanish Ministerio de Econom\'{\i}a y Competitividad---Fondo Europeo de Desarrollo Regional (FEDER) MTM 2010--15223 and by the Junta de Andaluc\'{\i}a grants FQM-336, FQM-1215 and FQM-246}
\author{Alberto Elduque${}^\dagger$}
 \address{Department of Mathematics, University of Zaragoza}
    \email{elduque@posta.unizar.es}
\thanks{${}^\dagger$ Supported by the Spanish Ministerio de Econom\'{\i}a y Competitividad---Fondo Europeo de Desarrollo Regional (FEDER) MTM2010-18370-C04-02 and by the Diputaci\'on General de Arag\'on---Fondo Social Europeo (Grupo de Investigaci\'on de \'Algebra)}

\begin{abstract} Some fine gradings on the exceptional Lie algebras $\mathfrak{e}_6$, $\mathfrak{e}_7$ and $\mathfrak{e}_8$ are described. This list tries to be exhaustive.

\medskip


\end{abstract}

\maketitle


\section*{Introduction}

Gradings by groups have played a key role in the study of Lie algebras, and contributed to understand their structural properties. Several families of examples can be found in \cite{librogradings}. To begin with, the root space decomposition of a complex semisimple Lie algebra is a grading by the group   $\mathbb{Z}^r$, with $r$ the rank of the Lie algebra. Any grading by a torsion-free abelian group is equivalent to a coarsening of such root space decomposition, and these gradings have been extensively used in representation theory. Gradings by not necessarily reduced root systems are very nice examples of this situation.
In particular, gradings by the integers have had frequent applications to physics, and they are specially relevant in algebraic contexts: if $J$ is a Jordan algebra, the Tits-Kantor-Koecher construction applied to $J$ is a
$\mathbb{Z}$-graded Lie algebra
 $L=L_{-1}\oplus L_{0}\oplus L_{1}$ with $L_{1}=J$, and the product in $J$ can be recovered from the one in $L$.
  Some other Jordan systems are related to `longer' $\mathbb{Z}$-gradings too.

Gradings by groups with torsion are also ubiquitous: gradings by cyclic groups and the corresponding finite order automorphisms are described by Kac \cite{Kacordenfinito}. They are intimately related to (infinite-dimensional) Kac-Moody Lie algebras and gradings by $ \mathbb{Z}$ on them. In Differential Geometry, they are  connected with symmetric spaces and their generalizations. Gradings by finite abelian groups are related to Lie color algebras (a generalization of Lie superalgebras), and sometimes to Lie algebra contractions. An interested reader can consult \cite[Introduction]{librogradings}  for these and other examples,  as well as for applications and references. Also a good and deep compilation of results can be found in \cite[Chapter~3, \S3]{enci}, which deals with the relationship between gradings on complex semisimple Lie algebras and automorphisms, and exhibit a wide variety of examples of gradings.

J.~Patera and H.~Zassenhaus, convinced about the relevance of gradings, initiated in \cite{LGI} a systematic study of  gradings on Lie algebras, 
emphasizing the role of the so called fine gradings (gradings which cannot be further refined).
Since then, a considerable number of authors have been trying to obtain a classification of  gradings on the simple Lie algebras (see, e.g. \cite{LGII,BahtKotclasicas,Albertoclasicas,g2,f4,e6}), which has culminated in the recent monograph \cite{librogradings}, where gradings on the classical simple Lie algebras and on the exceptional simple Lie algebras of types $G_2$ and $F_4$ are thoroughly studied. However, there is still work to be done. On one hand, not much is known about gradings on solvable or nilpotent Lie algebras. On the other hand,  the classification of gradings (for instance, the classification of fine gradings up to equivalence) is not yet finished for the complex exceptional simple Lie algebras of types $E_7$ and $E_8$ (denoted by $\mathfrak{e}_7$ and $\mathfrak{e_8}$), and this is also the case for the simple Lie algebras of types $D_4$ and $E_r$ ($r=6,7,8$) over algebraically closed fields of prime characteristic. Over the real numbers, even the classification of fine gradings for the classical Lie algebras is missing, although many low-dimensional cases have been considered.

Our goal is the classification of the fine gradings on the exceptional Lie algebras  $\mathfrak{e}_7$ and $\mathfrak{e}_8$  over an algebraically closed field $\FF$ of characteristic zero, and hence to finish the classification of fine gradings on simple Lie algebras over the complex numbers. (Note the result in \cite[Proposition~2]{e6}, which shows that the complex case yields a solution over arbitrary algebraically closed fields of characteristic $0$.)

This goal has not been reached yet. The purpose of this paper is to describe a  list of known fine gradings, which are compiled in our Main Theorem (Theorem~\ref{th:Main}) on the exceptional simple Lie algebras of type $E$.
Most of these gradings make sense in much more general contexts but, to avoid confusion, we will restrict ourselves to an algebraically closed ground field $\FF$ of characteristic zero.
The list exhausts the fine gradings, up to equivalence, on $\mathfrak{e}_6$, and we conjecture that it also exhausts them for $\mathfrak{e}_7$ and $\mathfrak{e}_8$. This has been announced in \cite[Figure~6.2]{librogradings}, although not all the descriptions there coincide with ours. Further details will be provided here too.

In our setting, a fine grading is the eigenspace decomposition relative to a maximal quasitorus (that is, a maximal abelian diagonalizable, or MAD, subgroup) of the group of automorphisms. This means that fine gradings are related to `large' abelian groups of symmetries. Hence our goal is equivalent to the classification of the MAD-subgroups, up to conjugation, of $\Aut(\mathfrak{e}_7)$ and $\Aut(\mathfrak{e_8})$.

Some of these MAD-subgroups have appeared in the literature.
For instance,  Griess describes in \cite{Griess} the  maximal  elementary $p$-subgroups of the groups $E_6$, $E_7$ and $E_8$. 
A larger class of abelian subgroups (not just MAD-subgroups) is studied in \cite{Yupr}.
The $\mathbb{Z}_6^3$-subgroup of $E_8$ has been considered by Hang and Vogan in \cite[pp.~22-25]{salelaZ6cubo}.
Also, Alekseevskii described the Jordan finite commutative subgroups of the groups $E_6$, $E_7$ and $E_8$ in \cite{Alek}. These include a MAD-subgroup of $E_8$ isomorphic to $\mathbb{Z}_5^3$. This subgroup has gained some attention lately, as it appeared in a talk by Kostant (see \cite{dondelaZ52deE8}) about the controversial paper by the  physicist Lisi \cite{Lisi}, which proposed a theory to go beyond the Standard Model in that it unifies all 4 forces of nature by using as gauge group the exceptional Lie group $E_8$. Kostant's talk, strictly mathematical,  dealt about an elaboration of the mathematics of $E_8$  in order to refute Lisi's Theory.
This is one more evidence of the fascination produced by the richness of $E_8$, and shows the relevance of understanding as much as possible about  this group and its tangent Lie algebra. Note that not even the finite abelian maximal groups are conveniently well known (consult the recent work \cite{e8finito}).
In terms of gradings, some of our descriptions  have appeared in \cite{gradingssymcomalg}, which uses gradings on composition algebras to construct some \emph{nice} gradings on exceptional algebras (e.g. a $\mathbb{Z}_3^5$-grading and a $\mathbb{Z}_2^8$-grading  on $\mathfrak{e}_8$).

This paper then gathers a lot of known material, and describes it in a homogeneous way. It is an expanded version of the talk presented by the first author in the conference \emph{Advances in Group Theory and Applications 2013}.

The paper is structured as follows.  First there is a section to recall the background: basic concepts
about gradings and their connection with groups of automorphisms, and some algebraic structures involved in the description of the exceptional Lie algebras: composition (and symmetric composition) algebras, Jordan algebras and its generalizations, and structurable algebras. Second, we present several models of the exceptional Lie algebras and constructions leading to them. After reviewing slightly how the exceptional Lie algebras emerged, we focus mainly on three constructions due to Elduque, Kantor and Steinberg, as they provide a convenient way to describe the fine gradings we are interested in. Finally, the third section describes some fine gradings on $\mathfrak{e}_6$, $\mathfrak{e}_7$ and $\mathfrak{e}_8$ starting from gradings on the `coordinate algebras' involved in the constructions above. (These algebras are usually much simpler than $\mathfrak{e}_6$, $\mathfrak{e}_7$ and $\mathfrak{e}_8$.) Fourteen fine gradings will be given on each of these simple Lie algebras.
The conjecture arises immediately of whether this is the complete list of fine gradings, up to equivalence, for these algebras.

\smallskip

For the sake of clarity of the exposition, in what follows the ground field $\FF$ will be assumed to be algebraically closed of characteristic $0$, even though many results are valid in more general contexts.


\section{Preliminaries}


\subsection{Gradings and automorphisms}

We begin by recalling the basics about gradings.
Let $ \mathcal{A}$ be a finite-dimensional algebra (not necessarily associative) over $\FF$, and let $G$ be an abelian group.

\begin{Definition}
A \emph{$G$-grading} $\Gamma$ on $ \mathcal{A}$ is a vector space decomposition
$$
 \Gamma:  \mathcal{A }= \bigoplus_{g\in G}  \mathcal{A}_g
$$
such that
$
\mathcal{A}_g  \mathcal{A}_h\subseteq  \mathcal{A}_{g+h}
$
for all $g,h\in G$.
\end{Definition}

Once such a decomposition is fixed, the algebra $ \mathcal{A}$ will be called a \emph{$G$-graded algebra}, the subspace $\mathcal{A}_g$ will be referred to as the \emph{homogeneous component of degree $g$} and its nonzero elements will be called the \emph{homogeneous elements of degree $g$}. The \emph{support} is the set $\supp \Gamma :=\{g\in G\mid \mathcal{A}_g\neq 0\}$.

\begin{Definition}
If $\Gamma\colon \mathcal{A}=\oplus_{g\in G} \mathcal{A}_g$ and $\Gamma'\colon \mathcal{A}=\oplus_{h\in H} \mathcal{A}_{h}$ are gradings by two abelian groups $G$ and $H$, $\Gamma$ is said to be a \emph{refinement} of $\Gamma'$ (or $\Gamma'$ a \emph{coarsening} of $\Gamma$) if for any  $g\in G$ there is $h\in H$ such that $\mathcal{A}_g\subseteq \mathcal{A}_{h}$. In other words, any homogeneous component of $\Gamma'$ is the direct sum of some homogeneous components of $\Gamma$. A refinement is \emph{proper} if some inclusion $\mathcal{A}_g\subseteq \mathcal{A}_{h}$ is proper. A grading is said to be \emph{fine} if it admits no proper refinement.
\end{Definition}

\begin{Definition}
Let $\Gamma$ be a $G$-grading on $\mathcal{A}$ and $\Gamma'$ an $H$-grading on another algebra $\mathcal{ B}$, with supports, respectively, $S$ and $T$. Then $\Gamma$ and $\Gamma'$ are said to be \emph{equivalent} if there is an algebra isomorphism  $\varphi\colon\mathcal{A} \rightarrow  \mathcal{B}$ and a bijection $\alpha\colon  S \rightarrow T$ such that $\varphi(\mathcal{A}_s)= \mathcal{B}_{\alpha(s)}$ for all $s\in S$. Any such $\varphi$ is called an \emph{equivalence} of $\Gamma$ and $\Gamma'$.
\end{Definition}

The study of gradings is based on classifying fine gradings up to equivalence, because any grading is obtained as a coarsening of some fine one.
 We will make use of the following invariant by equivalences:
 \begin{Definition}
 The \emph{type} of a grading $\Gamma$ is the sequence of numbers $(h_1,\ldots, h_r)$ where $h_i$ is the number of homogeneous components of dimension $i$, with $i=1,\ldots, r$ and $h_r\neq 0$. Obviously, $\dim \mathcal{A} = \sum_{i=1}^r ih_i$.
\end{Definition}

Given a group grading $\Gamma$ on an algebra $ \mathcal{A}$, there are many groups $G$ such that $\Gamma$, regarded as a decomposition into a direct sum of subspaces such that the product of any two of them lies in a third one, can be realized as a  $G$-grading, but there is one distinguished group among them (\cite{LGI}). Define $U(\Gamma)$ as the abelian group generated by $S=\supp\Gamma$ with defining relations $s_1s_2=s_3$ whenever $0\ne\mathcal{A}_{s_1}\mathcal{A}_{s_2}\subset\mathcal{A}_{s_3}$ ($s_i\in S$). It is called the
\emph{universal group of $\Gamma$}, since it verifies that, for any other realization of $\Gamma$ as a $G$-grading, there exists a unique homomorphism $U(\Gamma)\to G$ that restricts to identity on $\supp\Gamma$.
All the gradings throughout this work will be considered by their universal groups.
\smallskip

The classification of fine gradings on $\mathcal{A}$, up to equivalence, is the same as the classification of maximal diagonalizable subgroups (i.e., maximal quasitori) of $\Aut(\mathcal{A})$, up to conjugation (see e.g. \cite{LGI}). More precisely, given a $G$-grading  on the algebra $\mathcal{A}=\oplus_{g\in G} \mathcal{A}_g$, any $\chi$ belonging to
the group of characters $\hat G=\Hom(G,\FF^\times)$, acts as an automorphism of $\mathcal{A}$ by means of $\chi.x=\chi(g)x$ for any $g\in G$ and $x\in \mathcal{A}_g$. In case $G$ is the universal group of the grading, this allows us to identify $\hat G$ with a quasitorus (the direct product of a torus and a finite subgroup) of the algebraic group $\Aut(\mathcal{A})$. This quasitorus is the subgroup $\Diag(\Gamma)$ consisting of the automorphisms $\varphi$ of $\mathcal{A}$ such that the restriction of $\varphi$ to any homogeneous component is the multiplication by a (nonzero) scalar.
(See \cite[Chapter~3, \S 3]{enci} or \cite[\S 1.4]{librogradings}.)
Conversely, given a quasitorus $Q$ of $\Aut(\mathcal{A})$,
then $Q$ induces a  $\hat Q$-grading on $\mathcal{A}$,
where $\mathcal{A}_g=\{x\in \mathcal{A}\mid \chi(x)=g(\chi)x \,\forall\chi\in Q\}$ for any $g\in \hat Q$.  In this way   the fine gradings on $\mathcal{A}$, up to equivalence, correspond to the conjugacy classes in $\Aut(\mathcal{A})$ of the   maximal abelian diagonalizable subgroups  of $\Aut(\mathcal{A})$.


\subsection{Related structures}

We will recall here some algebraic structures involved in the constructions of the exceptional Lie algebras.
A very nice introduction to  nonassociative algebras  can be found in \cite{Schafernoasociativo}, but  the necessary material of composition algebras and Jordan algebras  is included here for completeness, as well as material about symmetric composition algebras and structurable algebras, which are not so well known.

\subsubsection{Composition algebras}\label{subsec_Compositionalgebras}

A Hurwitz algebra over $\FF$   is a  unital algebra $C$ endowed with a nonsingular quadratic form $q\colon C\to \FF$ which is multiplicative, that is, $q(xy)=q(x)q(y)$. This form $q$ is usually called the \emph{norm}.
Each element $a\in C$ satisfies
$$
a^2-t_C(a)a+q(a)1=0,
$$
where $t_C(a)=q(a+1)-q(a)-1$ is called the \emph{trace}. Denote $C_0=\{a\in C\mid t_C(a)=0\}$ the subspace of traceless elements. Note that $[a,b]=ab-ba\in C_0$ for any $a,b\in C$, since $t_C(ab)=t_C(ba)$.
The map $-\colon C\to C$ given by $\bar a=t_C(a)1-a$ is an involution and $q(a)=a\bar a$ holds.

There are Hurwitz algebras only in dimensions $1$, $2$, $4$ and $8$ (see e.g. \cite{librodelosrusos}).
 Moreover, under our hypothesis on the field,  there is only one (up to isomorphism) Hurwitz algebra of each possible dimension,
namely:
\begin{itemize}
\item the ground field $\FF$, with $q(a)=a^2$;
\item $\FF\times \FF$, with  componentwise product and norm given by $q(a,b)=ab$;
\item $\Mat_{2\times2}(\FF)$, with the usual matrix product and norm given by the determinant;
\item the split Cayley algebra over $\FF$.  This algebra can be characterized by the existence of a basis
$\{e_1,e_2,u_1,u_2,u_3,v_1,v_2,v_3\}$, which we call \emph{standard basis},  with multiplication given by
$$
\begin{array}{lll}
e_1u_j=u_j=u_je_2,\ &  u_iv_i=e_1, \    &    u_iu_{i+1}=v_{i+2}=-u_{i+1}u_i,       \\
e_2v_j=v_j=v_je_1,&   v_iu_i=e_2,   &     -v_iv_{i+1}=u_{i+2}=v_{i+1}v_i,
\end{array}
$$
all the remaining products being $0$,
and the   polar form of the norm (also denoted by $q$) of two basic elements is zero except for $q(e_1,e_2)=1=q(u_i,v_i)$, $i=1,2,3$.
\end{itemize}
With the exception of the Cayley algebra, all of these are associative. The Cayley algebra is not associative but alternative (the algebra generated by any pair of elements is associative). We will use the notations $\FF$, $\K$, $\Q$ (usually called \emph{quaternion algebra}) and $\C$ (usually called \emph{octonion algebra}) for each of these algebras, respectively.

Recall that for any $a,b\in C$, the endomorphism
$$
d_{a,b}=[l_a,l_b]+[l_a,r_b]+[r_a,r_b]
$$
is a derivation of $C$ for $l_a(b)=ab$ and $r_a(b)=ba$. This will be instrumental to construct Lie algebras from Hurwitz algebras.

\subsubsection{Jordan algebras}\label{subsec_Jordanalgebras}

A \emph{Jordan algebra} is a commutative (nonassociative) algebra satisfying the Jordan identity
$$
(x^2y)x=x^2(yx).
$$
This kind of algebras were introduced by Jordan in 1933 to formalize the notion of an algebra of observables in quantum mechanics. Such line of research was abandoned time ago, but Jordan algebras have found a range of applications because of their relationship to Lie algebras. A standard reference is \cite{JacobsondeJordan}.

If $A$ is an associative algebra (with multiplication denoted by juxtaposition) and we consider the new product on $A$ given by
$$
x\circ y=\frac12(xy+yx),
$$
then $(A,\circ)$ is a Jordan algebra, denoted by $A^+$. A Jordan algebra which is a subalgebra of $A^+$ for some associative algebra $A$, is called a \emph{special} Jordan algebra, and otherwise it is called \emph{exceptional}.
If $(A,-)$ is an associative algebra with involution, then the set of hermitian elements $H(A,-)=\{a\in A\mid \bar a=a\}$ is a subalgebra of $A^+$ (not of $A$), and hence it is a special Jordan algebra.

In particular, if $C$ is an associative Hurwitz algebra with involution given by $-$, the algebra
  $ H_n(C,*)=\{x=(x_{ij})\in \Mat_{n\times n}(C)\mid x_{ij}=\bar x_{ji}\}$ is a Jordan algebra for any $n\ge 3$. (For $n=1$ or $n=2$ this is also true, but in a trivial way, so we will assume $n\geq 3$.) It is proved in \cite{JacobsondeJordan} that if $\C$ is the Cayley algebra, $H_n(\C,*)$ is a Jordan algebra if and only if $n=3$. Besides, this is the only exceptional Jordan algebra, which is called the Albert algebra, and  will be denoted  by $\mathbb{A}$.

  If $J=H_n(C,*)$ for some Hurwitz algebra and some $n$, consider the linear map $t_J\colon J\to \FF$ given by $t_J(x)=\frac{\tr(x)}{n}=\frac{\sum_{i=1}^n x_{ii}}{n}$. This map is called the \emph{normalized trace} and it is the only linear map such that $t_J(1)=1$ and $t_J((xy)z)=t_J(x(yz))$ for any $x,y,z\in J$. Thus we have a decomposition $J=\FF1\oplus J_0$, for $J_0=\{x\in J\mid t_J(x)=0\}$, since $x*y=xy-t_J(xy)1\in J_0$. In particular we have a commutative multiplication $*$ defined in $J_0$.

  If $J$ is a Jordan algebra and $R_x\colon J\to J$, $y\mapsto yx$ is the multiplication operator, observe that
  \begin{equation}\label{eq_productosdeRsenJordan}
  [[R_x,R_y],R_z]=R_{(yz)x-y(zx)}
  \end{equation}
   for any $x,y,z\in J$,
  and thus, the \emph{structure algebra} $\mathfrak{Str}(J)$, or Lie algebra generated by the multiplication operators, coincides with
  $R_J+[R_J,R_J]$ (this sum is direct if $J$ is unital). It is also a consequence of Equation~(\ref{eq_productosdeRsenJordan}) that $[R_J,R_J]$ is an ideal of the Lie algebra of derivations
  $\Der(J)=\{d\in\mathfrak{gl}(J)\mid d(xy)=d(x)y+xd(y)\ \forall x,y\in J\}$. The algebra generated by the traceless multiplication operators $\{R_x\mid x\in J_0\}$ is called the \emph{inner structure algebra} and it also coincides with
  $R_{J_0}+[R_J,R_J]$.

  In case $J=\mathbb{A}$ is the Albert algebra,
  $\Der(\mathbb{A})$ is simple \cite[Theorem~3]{Jacobsondeexcepcionales}, so in particular every derivation is inner ($\Der(\mathbb{A})=[R_{\mathbb{A}},R_{\mathbb{A}}]$). The inner structure algebra is simple too \cite[Theorem~4]{Jacobsondeexcepcionales}, and these provide our first two models of $\mathfrak{f}_4$ and $\mathfrak{e}_6$. Moreover, any element $x\in\mathbb{A}$ satisfies a cubic equation
  \begin{equation}\label{eq_laecuacioncubicaAlbert}
  x^3-T(x)x^2+Q(x)x-N(x)1=0,
  \end{equation}
  for the scalars $T(x)=\tr(x)$, $Q(x)=\frac12((T(x))^2-T(x^2))$
  and $N(x)=\frac16((T(x))^3-3T(x)T(x^2)+2T(x^3))$.
  The cubic form $N$ is also closely related to $\mathfrak{e}_6$.

\subsubsection{Symmetric composition algebras}\label{subsec_Symmetriccompositionlgebras}

A  \emph{symmetric composition algebra} is a triple $(S,*,q)$, where $(S,*)$ is a (nonassociative) algebra over $\FF$ with multiplication denoted by $x*y$ for $x,y\in S$, and where $q\colon S\to \FF$ is a regular quadratic form verifying
\begin{equation*}
\begin{array}{ll}
q(x*y)=q(x)q(y),\\
q(x*y,z)=q(x,y*z),
\end{array}
\end{equation*}
for any $x,y,z\in S$, where $q(x,y)=q(x+y)-q(x)-q(y)$ is the polar form of $q$.

\begin{Example}
Let $C$ be a Hurwitz algebra with norm $q$. The same vector space with new product
$$
x*y=\bar x\bar y
$$
for any $x,y\in C$ is a symmetric composition algebra for the same norm, called the \emph{para-Hurwitz algebra} attached to the Hurwitz algebra $C$. We will denote it by $pC=(C,*,q)$. 
Note that the unit of $C$ becomes
a \emph{paraunit} in $pC$, that is, an element $e$ such that $e * x = x * e = q(e, x)e - x$.
\end{Example}

\begin{Example}
The \emph{Okubo algebra},  or pseudo-octonion algebra, is the algebra
defined on the subspace of trace $0$ matrices of degree $3$: $\Ok=(\Mat_{3\times3}(\FF)_0,\ast,q)$
with multiplication
\begin{equation}\label{eq_productoOkubo}
x*y=\omega xy-\omega^2 yx-\frac{\omega-\omega^2}{3}\tr(xy)1
\end{equation}
and norm $q(x):=\frac16\tr(x^2)$, for $x,y\in \Mat_{3\times3}(\FF)_0$, where $\omega$ is a primitive cubic root of $1$.
This algebra is a symmetric composition algebra, but   it does not have an identity element (and it is not alternative). It was introduced by Okubo in \cite{Okubodef}, who was working in Particle Physics and the symmetry
given by the compact group  $\textrm{SU}(3)$
(the real algebra $\mathfrak{su}(3) = \{x\in \mathfrak{sl}(3,\mathbb{C}) \mid x^* = -x\}$, for $*$ the unitary involution, is closed for the product given in Equation~(\ref{eq_productoOkubo})).
\end{Example}

The classification of the symmetric composition algebras was obtained in \cite{AlbSymmCompalgebras}.

\begin{Theorem}
Every symmetric composition algebra over $\FF$ (algebraically closed) is isomorphic either to a para-Hurwitz algebra or to the Okubo algebra. That is, there are only five symmetric composition algebras (up to isomorphism), namely, $\pF$, $\pK$, $\pQ$, $\pC$ and $\Ok$.
\end{Theorem}

\subsubsection{Structurable algebras}\label{subsec_Structurablealgebras}

Let $(A,-)$ be a unital algebra with involution \lq\lq$-$\rq\rq. Denote the multiplication in $A$ by juxtaposition.
For $x,y\in A$, consider the linear operator $V_{x,y}\in \End(A)$ given by $V_{x,y}(z)=(x\bar y)z+(z\bar y)x-(z\bar x)y\equiv \{x,y,z\}$. The algebra $A$ is called a \emph{structurable algebra} in case the identity
$$
\{x,y,\{z,w,v\}\}=\{\{x,y,z\},w,v\}-\{z,\{y,x,w\},v\}+\{z,w,\{x,y,v\}\}
$$
is satisfied for any $x,y,z,w,v\in A$, or equivalently,
$$
{[}V_{x,y},V_{z,w}]=V_{V_{x,y}z,w}-V_{z,V_{y,x}w}.
$$
The reader may consult \cite{Allisonestructurables} for the definition and properties of structurable algebras.

\begin{Example}
Any (unital) associative algebra  with involution $(A,-)$ is a structurable algebra.
\end{Example}

The space $\mathfrak{Instr}(A,-)=\{\sum_i V_{x_i,y_i}\mid x_i,y_i\in A\}$
is a subalgebra of the Lie algebra $\mathfrak{gl}(A)$, called the \emph{inner structure algebra} of $(A,-)$.
The map
$$
\epsilon\colon \mathfrak{Instr}(A,-)\to \mathfrak{Instr}(A,-),
\qquad\epsilon(V_{x,y})=-V_{y,x}
$$
 if $x,y\in A$, is an involutive automorphism of this Lie algebra.
 Thus $\mathfrak{Instr}(A,-)$ turns out to be $\mathbb{Z}_2$-graded. 
The elements in   $H(A,-)=\{x\in A\mid \bar x=x\}$ and $S(A,-)=\{x\in A\mid \bar x=-x\}$ are called hermitian
and skew-hermitian respectively. It follows that $\mathfrak{Instr}(A,-)_{\bar 0}=\bigl(\mathfrak{Instr}(A,-)\cap \Der(A,-)\bigr)\oplus V_{S(A,-),1}$ and $\mathfrak{Instr}(A,-)_{\bar 1}=V_{H(A,-),1}$,
where $\Der(A,-)$ denotes the set of derivations of $A$ that commutes with the involution$-$.

\begin{Example}\label{ex_estructurableJordan}
If $J$ is a Jordan algebra, then $(J,-)$ is a structurable algebra with the involution $-$ given by the identity map. In this case  $V_{x,y}=R_{x y}+[R_x,R_y]$ for $x,y\in J$, where $R_x$ is the   multiplication operator by $x$. In this sense, the inner structure algebra of $(J,-)$ is the usual inner structure algebra for a Jordan algebra described in Subsection~\ref{subsec_Jordanalgebras}, and the  $\mathbb{Z}_2$-grading produced by $\epsilon$ is $\mathfrak{Instr}(J,-)_{\bar0}=[R_J,R_J]$, $\mathfrak{Instr}(J,-)_{\bar1}= R_{J_0}$.
\end{Example}

\begin{Example}\label{ex_tensoresalgcomposicion}
If $(C_1,-)$ and $(C_2,-)$ are composition algebras over the field $\FF$, then $(C_1\otimes C_2,-)$ is a structurable algebra (see \cite[Example~6.6]{paraSteinberg})
for the product given by $(a\otimes b)(c\otimes d)=ac\otimes bd$ and the involution given by
$\overline{a\otimes b}=\bar a\otimes \bar b$.
\end{Example}


\section{Models of the exceptional Lie algebras}

The study of the Lie algebras began at the end of the XIX century, once Lie had translated certain problem of transformation groups to an algebraic context. The first fundamental contributions are due to Killing, who classified the complex simple Lie algebras in four key papers published during the years 1886 and 1890. Initially he thought that the only possible cases were the Lie algebras of the special linear group $\textrm{SL}_n(\mathbb{C})$ and of the orthogonal and symplectic groups $\textrm{O}(n,\mathbb{C})$ and $\textrm{Sp}(n,\mathbb{C})$, now called the classical Lie algebras. But, during his work, he obtained that, besides the classical Lie algebras, there were a few other Lie algebras, of dimensions $78$, $133$, $248$, $52$ and $14$, now denoted as $\mathfrak{e}_6$, $\mathfrak{e}_7$, $\mathfrak{e}_8$, $\mathfrak{f}_4$ and $\mathfrak{g}_2$. Actually, he only proved the existence for  $\mathfrak{g}_2$, but he described all possibilities for rank, dimension and root systems. He found six algebras, since he did not notice that two of them were isomorphic (case $\mathfrak{f}_4$). This is a marvelous result, but the standard reference for it   is Cartan's thesis in 1894, which completed Killing's classification,
giving a rigorous treatment.
This is a fundamental contribution, where Cartan proved the existence of all the exceptional simple Lie algebras.

\subsection{First Models}

The history of these algebras has been growing in parallel to the one of the related Lie groups.
The first description of the smallest of the exceptional Lie groups was due to Engel (\cite{EngelG2}), who, in 1900,
described it  as the isotropy
group of a generic 3-form in 7 dimensions. 
\'{E}lie Cartan was the first to consider the group $G_2$   as the automorphism
group of the octonion algebra in 1914 \cite[p.~298]{Cartang2} (although he commented about it earlier),
as well the Lie algebra $\mathfrak{g}_2$ as the derivation algebra of the octonions (both on the split and compact forms). Jacobson generalized this result to arbitrary fields (\cite{JacobsonG2}).
  This approach
  became popular
through the article \cite{Freudenthal} by Hans Freudenthal,    in 1951.
But, for rather a long time, $G_2$ was the only Lie
group for which further results were obtained.

The following  model of a exceptional Lie algebra had to wait until 1950,
when Chevalley and
  Schafer (\cite{ChevalleySchafer})  showed that  the set of derivations of the Albert algebra $\mathbb{A}$ is $\mathfrak{f}_4$. Tomber proved in \cite{Tomber} 
  the converse: a Lie algebra over a field of characteristic $0$ is of type $F_4$  if and only if it is isomorphic to the derivation algebra of an exceptional
  simple Jordan algebra.
  This fact led Tits, among other authors, to study the relationship between
  Jordan algebras and the remaining exceptional simple Lie algebras, which were constructed in
  a unified way \cite{Tits}. We will revise this construction in the following subsection.

 The algebra $ \mathfrak{e}_6$ is also closely related to the Albert algebra. On one hand, it is the inner structure algebra of the Albert algebra  (the Lie algebra generated by the right multiplication operators $R_a$ for $a\in \mathbb{A}$ with zero trace, as in Subsection~\ref{subsec_Jordanalgebras}).
 On the other hand, if $N\colon\mathbb{A}\to \FF$ is the cubic norm as in Equation~(\ref{eq_laecuacioncubicaAlbert}), and  $N(a,b,c)$ denotes the trilinear form obtained by polarization, then the Lie algebra  $ \mathfrak{e}_6$ can be characterized as the Lie algebra $\{f\in \End(\mathbb{A})\mid N(f(a),b,c)+N(a,f(b),c)+N(a,b,f(c))=0\}$ (\cite{JacobsonE6} and \cite{Freudenthal}).

 The first model of $ \mathfrak{e}_7$  related to the Albert algebra was provided in  \cite{TitsE7}, as the algebra defined on the vector space $\left(\mathbb{A}\otimes \mathfrak{sl}_2(\FF)\right) \oplus\Der(\mathbb{A})$ with the product
 $$\begin{array}{rl}
 [x\otimes a+d_1,y\otimes b+d_2]=&xy\otimes [a,b]+d_2(x)\otimes a-d_1(y)\otimes b\\
 &-\frac12\tr(\ad a\ad b)[R_x,R_y]+[d_1,d_2],\end{array}
 $$
 for $x,y\in \mathbb{A}$, $a,b\in\mathfrak{sl}_2(\FF)$ and $d_1,d_2\in\Der(\mathbb{A})$.
 The details of this construction appear in \cite[\S9]{Jacobsondeexcepcionales}.
 This is a version of  what nowadays is called the Tits-Kantor-Koecher construction applied to the Albert algebra. The name refers to several constructions which appeared almost simultaneously, and turned out to be essentially equivalent. In Koecher's construction \cite{Koecher}, one forms $\mathbb{A}\oplus\bar{\mathbb{A}}\oplus \mathfrak{Str}(\mathbb{A})$, where $\bar{\mathbb{A}}$ is simply a copy of the vector space $\mathbb{A}$, with the anticommutative product given by $[x,y]=0=[\bar x,\bar y]$, $[x,\bar y]=2R_{xy}+2[R_y,R_x]$ if $x,y\in \mathbb{A}$, and $[L,x]=L(x)$, $[L,\bar x]=\overline{\bar L(x)}$, if $L=R_x+\sum_i[R_{x_i},R_{y_i}]\in \mathfrak{Str}(\mathbb{A})$, where $\bar L=-R_x+\sum_i[R_{x_i},R_{y_i}]$. This construction will be generalized in Subsection~\ref{subsec_Kantorsmodels}.
Similar to the situation for $\mathfrak{e}_6$, the Lie algebra  $\mathfrak{e}_7$  can be characterized too as the set of linear transformations of certain vector space $M$ leaving invariant a quartic form \cite{FreudenthalE7}. Here, as a vector space, $M$ is $ \mathbb{A}\oplus\bar{\mathbb{A}}\oplus \FF\oplus\bar{\FF}$. This will play an important role in our description of the gradings with automorphisms of order $4$ involved, in Subsections~\ref{subsec_Z43enlaestructurable} and \ref{subsec_gradingsconcuatros}.

Finally, the difficulty of finding a good model for  $\mathfrak{e}_8 $ (coordinate free, that is, not given by means of generators and relations obtained from the root system) is that the nontrivial  representation  of minimal dimension for $\mathfrak{e}_8$ is the adjoint representation, so there is no hope to embed  $\mathfrak{e}_8 $  as a subalgebra of $\mathfrak{gl}(V)$ for some vector space $V$ of smaller dimension. However, some other  `linear models' can help in this purpose. Let $V$ be a vector space of dimension $9$, then we can construct  $\mathfrak{e}_8 $ as the vector space
$$
\bigwedge^3V^*\oplus \mathfrak{sl}(V)\oplus\bigwedge^3V,
$$
with Lie bracket as in \cite[Exercise\,22.21]{FultonHarris}, based on the trilinear map given by the usual wedge product $\bigwedge^3V\otimes\bigwedge^3V\otimes\bigwedge^3V\to\bigwedge^9V\cong\FF$.
\smallskip

We refer to \cite{Jacobsondeexcepcionales,Schafernoasociativo,Freudenthalcuadradomagico,Freudenthal} for these and other algebraic constructions of the exceptional Lie algebras. We stress the reference \cite{enci}, where many models appear:
\cite[Chapter~5, \S1]{enci} is devoted to models of exceptional Lie algebras associated to a Cayley algebra (over arbitrary fields of characteristic zero, with several references to the reals), while \cite[Chapter~5, \S2]{enci}
provides other models based on gradings.
This shows the interactions between nice models and certain gradings, and this philosophy is certainly present in Subsections~\ref{subsec_elprimo5} and \ref{subsec_explicacioneslinealesdeloscuatros}.


\subsection{Tits construction}

In 1966, Tits gave a unified construction of the exceptional simple Lie algebras (over fields of characteristic not two and three) in \cite{Tits}. The construction used   a couple of ingredients: an alternative algebra of degree 2 and a Jordan algebra of degree 3. In case these ingredients are chosen to be
Hurwitz algebras    and   Jordan algebras of hermitian $3\times3$ matrices over Hurwitz algebras, Freudenthal's magic square \cite{Freudenthalcuadradomagico} is obtained. We recall the construction in our concrete case.

Let $C$ be a Hurwitz algebra over $\FF$ with norm $q$, and let  $J=H_3(C',*)$ be the Jordan algebra of hermitian $3\times 3$-matrices over another Hurwitz algebra $C'$.
Consider the vector space
$$
\mathcal{T}(C,J)=\Der(C)\oplus (C_0\otimes J_0)\oplus \Der(J)
$$
with anticommutative multiplication specified by
\begin{itemize}
\item $\Der(C)$ and $\Der(J)$ are Lie subalgebras,
\item $[\Der(C),\Der(J)]=0$,
\item $[d,a\otimes x]=d(a)\otimes x$, $[D,a\otimes x]=a\otimes D(x)$,
\item $[a\otimes x,b\otimes y]=t_J(xy)d_{a,b}+[a,b]\otimes x*y+2t_C(ab)[R_x,R_y]$
\end{itemize}
for all $d\in\Der(C)$, $D\in\Der(J)$, $a,b\in C_0$ and $x,y\in J_0$, with the notations of Subsections~\ref{subsec_Compositionalgebras} and \ref{subsec_Jordanalgebras}.
 Now, using all the possibilities for $C$ and $C'$, we obtain
Freudenthal's Magic Square as follows (\cite{Tits}) (note that we have added a column with $J=\FF$ to obtain $G_2$ with the same construction):
 \begin{center}{
 \begin{tabular}{cc}
 &\null\qquad$J$\\[4pt]
 \lower 6pt\hbox{$C$}&
\begin{tabular}{c|ccccc|}
$\mathcal{T}(C,J)$&$\FF$&$H_3(\FF)$ & $H_3(\FF\times \FF)$ & $H_3(\Mat_{2\times 2}(\FF))$ &$H_3(\C)$  \cr
 \hline
 $\FF$&0&$ \mathfrak{a}_1$& $\mathfrak{a}_2$ & $\mathfrak{c}_3$ &$\mathfrak{f}_4$ \cr
 $\FF\times \FF$&0&$\mathfrak{a}_2$& $\mathfrak{a}_2\oplus \mathfrak{a}_2$ & $\mathfrak{a}_5$ &$\mathfrak{e}_6$ \cr
 $\Mat_{2\times 2}(\FF)$&$\mathfrak{a}_1$&$\mathfrak{c}_3$& $\mathfrak{a}_5$ & $\mathfrak{d}_6$ &$\mathfrak{e}_7$ \cr
 $\C$&$\mathfrak{g}_2$&$\mathfrak{f}_4$& $\mathfrak{e}_6$ & $\mathfrak{e}_7$ &$\mathfrak{e}_8$ \cr
 \hline
  \end{tabular}\end{tabular}}
  \end{center}

\smallskip

\subsection{Some symmetric constructions}

In spite of the apparent asymmetry in the use of the two Hurwitz algebras in Tits construction, the Magic Square is symmetric. This lead several authors to look for more symmetric constructions. Such an approach was taken by Vinberg in \cite[p.~177]{enci}, and interpreted by Barton and Sudbery in \cite{BartonSudbery} as a construction depending on two composition algebras and closely related to the triality principle. A similar construction was provided by Landsberg and Manivel in \cite{LandsbergManivel}, inspired by previous work of Allison and Faulkner   \cite{paraSteinberg}.

The construction we are going to recall here (and use later on), is the construction in \cite{cuadradomagicoAlb},
based on two symmetric composition algebras, which has turned to be very useful in finding fine gradings on exceptional Lie algebras  \cite{gradingssymcomalg}.

Let $(S,\ast, q)$ be a symmetric composition algebra and let
$$
\mathfrak{o}(S,q)=\{d\in \End_\FF(S)\mid q(d(x),y)+q(x,d(y))=0\,\forall x,y\in S\}
$$
be the corresponding orthogonal Lie algebra. Consider the subalgebra of $\mathfrak{o}(S,q)^3$ defined by
$$
\mathfrak{tri}(S,\ast,q)=\{(d_0,d_1,d_2)\in \mathfrak{o}(S,q)^3 \mid
d_0(x\ast y)=d_1(x)\ast y + x\ast d_2(y)\ \forall x,y\in S\},
$$
 which is called the \emph{triality Lie algebra}.
The order three automorphism $\vartheta$ given by
$$
\vartheta\colon  \mathfrak{tri}(S,\ast,q)\longrightarrow \mathfrak{tri}(S,\ast,q), \quad (d_0, d_1, d_2)\longmapsto (d_2, d_0, d_1),
$$
 is called the \emph{triality automorphism}. Take the element of $\mathfrak{tri}(S,\ast,q)$ (denoted by $\mathfrak{tri}(S)$ when it is no ambiguity) given by
$$
t_{x,y}:=\left(\sigma_{x,y},\frac{1}{2}q(x,y)id-r_x l_y,\frac{1}{2}q(x,y)id-l_x r_y\right),
$$
where $\sigma_{x,y}(z)=q(x,z)y-q(y,z)x$, $r_x(z)=z\ast x$, and $l_x(z)=x\ast z$ for any $x,y,z\in S$.

Let $(S,\ast,q)$ and $(S',\star ,q')$ be two symmetric composition algebras over $\FF$.
Consider the following vector space, which depends symmetrically on $S$ and $S'$:
$$
\mathfrak{g}(S,S') := \mathfrak{tri}(S,\ast, q)\oplus \mathfrak{tri}(S',\star, q')
\oplus (\bigoplus_{i=0}^2 \iota_i(S\otimes S'))
$$
where $\iota_i(S\otimes S')$ is just a copy of $S\otimes S'$ ($i=0,1,2$), and the anticommutative product
on $\mathfrak{g}(S,S')$ is determined by the following conditions:
\begin{itemize}
\item   $\mathfrak{tri}(S,\ast, q)\oplus \mathfrak{tri}(S',\star, q')  \text{ is a Lie subalgebra of }  \mathfrak{ g}(S,S')$;
\item $[(d_0,d_1,d_2), \iota_i(x\otimes x')]=\iota_i(d_i(x)\otimes x')$,
   $[(d'_0,d'_1,d'_2), \iota_i(x\otimes x')]=\iota_i(x\otimes d'_i(x'))$,  for any $(d_0,d_1, d_2)\in \mathfrak{tri}(S)$, $(d'_0,d'_1, d'_2)\in \mathfrak{tri}(S')$, $x\in S$ and $x'\in S'$;
\item    $[\iota_i(x\otimes x'), \iota_{i+1}(y\otimes y')]=\iota_{i+2}((x*y)\otimes(x'\star y'))$ (indices modulo 3), for any $i=0,1,2$,  $x,y\in S$ and $x',y'\in S'$;
\item   $[\iota_i(x\otimes x'),\iota_i(y\otimes y')] = q'(x',y')\vartheta^i(t_{x,y}) + q(x,y)\vartheta'^i(t'_{x',y'})\in\mathfrak{tri}(S)\oplus\mathfrak{tri}(S')$,
    for any $i=0,1,2$,  $x,y\in S$ and $x',y'\in S'$,
    $\vartheta$ and $\vartheta'$ being the corresponding triality automorphisms.
\end{itemize}
  The anticommutative algebra $\mathfrak{g}(S,S')$ defined in this way turns out to be  a Lie algebra (\cite[Theorem~3.1]{cuadradomagicoAlb}), and we recover  Freudenthal's Magic Square if symmetric composition algebras of all possible dimensions are considered:
 \begin{center}{
 \begin{tabular}{cc}
 &$\dim S$\\[2pt]
 \lower 5pt\hbox{$\dim S'$}&
\begin{tabular}{c|cccc|}
&1 & 2 & 4 &8  \cr
 \hline
 1&$\mathfrak{a}_1$& $\mathfrak{a}_2$ & $\mathfrak{c}_3$ &$\mathfrak{f}_4$ \cr
 2&$\mathfrak{a}_2$& $\mathfrak{a}_2\oplus \mathfrak{a}_2$ & $\mathfrak{a}_5$ &$\mathfrak{e}_6$ \cr
 4&$\mathfrak{c}_3$& $\mathfrak{a}_5$ & $\mathfrak{d}_6$ &$\mathfrak{e}_7$ \cr
 8&$ \mathfrak{f}_4$& $ \mathfrak{e}_6$ & $ \mathfrak{e}_7$ &$ \mathfrak{e}_8$ \cr
 \hline
  \end{tabular}\end{tabular}}
  \end{center}

If $C_1$ and $C_2$ are two Hurwitz algebras over $\FF$, the Lie algebra $\mathfrak{g}(pC_1,pC_2)$ is isomorphic to $\mathcal{T}(C_1,H_3(C_2,*))$ \cite[Theorem~6.25]{librogradings}.


\subsection{Kantor's construction}\label{subsec_Kantorsmodels}

Let $(A,-)$ be a structurable algebra. Denote by  $S\equiv S(A,-)$ its set of skew-hermitian elements.
Endow the $\mathbb{Z}$-graded vector space $\mathcal{K}=\mathcal{K}_{-2}\oplus\mathcal{K}_{-1}\oplus\mathcal{K}_{0}\oplus\mathcal{K}_{1}\oplus\mathcal{K}_{2}$,
with
$$\begin{array}{ll}
\mathcal{K}_{2}=S,\quad&\mathcal{K}_{1}=A,
\\
\mathcal{K}_{-2}=S\tilde{\ }&\mathcal{K}_{-1}=A\tilde{\ },
\end{array}\qquad \mathcal{K}_{0}=\mathfrak{Instr}(A,-),
$$
where $S{\tilde{\ }}$ and $A{\tilde{\ }}$ are simply copies of $S$ and $A$ respectively,
with a graded Lie algebra structure given by the  anticommutative multiplication such that $\mathfrak{Instr}(A,-)=V_{A,A}$ is a subalgebra and the following conditions hold:
\begin{equation}\label{eq_productoenKantor}
\begin{array}{l}
\begin{array}{ll}
{[}T,a]=T(a), &
{[}T,a\tilde{\ }]=(T^\epsilon a)\tilde{\ } ,\\
{[}T,s]= T(s)+s\overline{T(1)},\qquad&
{[}T,s\tilde{\ }]=(T^\epsilon (s)+s\overline{T^\epsilon (1)})\tilde{\ } ,
\end{array}\\
\begin{array}{l}
{[}a+s,a'+s']=2(s\bar s'-s'\bar s)\in \mathcal{K}_{2},\\
{[}a\tilde{\ }+s\tilde{\ },a'\tilde{\ }+s'\tilde{\ }]=2(s\bar s'-s'\bar s)\tilde{\ }\in \mathcal{K}_{-2},\\
{[}a +s ,a'\tilde{\ }+s'\tilde{\ }]=(-s'a)\tilde{\ }+L_sL_{s'}+2V_{a,a'}+(sa')\in\mathcal{K}_{-1}\oplus\mathcal{K}_{0}\oplus\mathcal{K}_{1},
\end{array}\end{array}
\end{equation}
for $T\in \mathfrak{Instr}(A,-)$, $a,a'\in A$, $s,s'\in S$,
where $L_s\colon A\to A$ denotes the left multiplication by $s\in S$
(so that $2L_sL_{s'}=V_{ss',1}-V_{s,s'}\in \mathfrak{Instr}(A,-)$).
Following  \cite[\S\,6.4]{paraTKK}, this  is a $\mathbb{Z}$-graded Lie algebra denoted by $\mathfrak{Kan}(A,-)$, which be called the \emph{Kantor construction} attached to the structurable algebra $A$.
The construction takes its name from Kantor, who introduced it first in \cite{Kantor}, although in a somewhat different way.

In connection with Example~\ref{ex_estructurableJordan}, if $J$ is a Jordan algebra, the Kantor construction
$\mathfrak{Kan}(J,-)$ coincides (up to isomorphism) with the classical Tits-Kantor-Koecher Lie algebra constructed from the Jordan algebra $J$.  

A necessary and sufficient condition for a Lie algebra to be isomorphic to the Kantor's construction attached to a structurable algebra is given by the existence of  a $\mathfrak{sl}_2$-triple  $\{e,h,f\}$  (that is, $[h,e]=2e,[h,f]=-2f,[e,f]=h$) in $L$
such that
$L$ is the direct sum of irreducible modules for $\langle\{e,h,f\}\rangle$ of dimensions 1, 3, and 5;
the only ideal of $L$ which centralizes $ \{e,h,f\} $ is $\{0\}$, and
$L$ is generated by the eigenspaces $2$ and $-2$ for $\ad h$ \cite[Theorem~6.10]{paraTKK}.

The relevance for our purposes comes from the fact that   exceptional Lie algebras are obtained in terms of Kantor's constructions attached to certain structurable algebras. Some examples of this situation are shown next:

\begin{Example}
Consider the tensor products of composition algebras: $\C=\FF\otimes\C$, $\K\otimes\C$, $\Q\otimes\C$, and $\C\otimes\C$. These are structurable algebras according to Example~\ref{ex_tensoresalgcomposicion}. Kantor's construction gives the following Lie algebras:
$$
\mathfrak{Kan}(\C)\cong\mathfrak{f}_4,\quad
\mathfrak{Kan}(\K\otimes\C)\cong\mathfrak{e}_6,\quad
\mathfrak{Kan}(\Q\otimes\C)\cong\mathfrak{e}_7,\quad
\mathfrak{Kan}(\C\otimes\C)\cong\mathfrak{e}_8,
$$
as stated in \cite{Kantor} (see, alternatively, \cite[\S8(c)]{Allisonparaderystr}).
\end{Example}

\begin{Example}\label{ex_CDprocesoparaestructurables}
Let $J=H_4(C,*)$ be the Jordan algebra defined in Subsection~\ref{subsec_Jordanalgebras} for $C$ a Hurwitz associative algebra. This algebra is a finite-dimensional simple Jordan algebra, with generic trace $t_J(x)=\frac14\tr(x)$. We can use a sort of Cayley-Dickson process  with $J$ in order to construct a structurable algebra as follows. Take a nonzero $\mu\in \FF$ and the algebra defined on the vector space
$$
A=J\oplus vJ,
$$
where $v$ is simply a convenient mark to indicate the cartesian product of two copies of $J$, with multiplication and involution given by
\begin{equation}\label{eq_CDestructurable}
\begin{array}{l}
(x_1+vx_2)(x_3+vx_4):=x_1x_3+\mu(x_2x_4^\theta)^\theta+v(x_1^\theta x_4+(x_2^\theta x_3^\theta)^\theta),\\[4pt]
\overline{x_1+vx_2}:=x_1-vx_2^\theta
\end{array}
\end{equation}
for any $x_i\in J$ for $i=1,2,3,4$,
where the multiplication in $J$ is denoted by juxtaposition, with $x^\theta:=-x+2t_J(x)1$ (so that $1^\theta=1$) for any $x\in J$.
The resulting algebra  $(A,-)$ is denoted by ${\CD}(H_4(C))$.
According to  \cite{CayleyDicksonparaestructurables} (see, alternatively, \cite[Example~6.7]{paraTKK}), this algebra   is structurable with
the extra property that the space of skew-hermitian elements has dimension $1$. The Lie algebra obtained by means of  Kantor's construction $\mathfrak{Kan}({\CD}(H_4(C)))$ is $\mathfrak{e}_6$, $\mathfrak{e}_7$ and $\mathfrak{e}_8$, respectively, in case   $\dim C=1,2$ or $4$.
\end{Example}


\subsection{Steinberg's construction}

Consider the following example \cite{paraSteinberg}.

\begin{Example}
Let  $(A,-)$ be a unital associative algebra with involution  and consider the unitary Lie algebra
$\mathfrak{u}_n(A,-)=\{x\in\Mat_{n\times n}(A)\mid \bar x^t=-x\}$. Some remarkable elements are $u_{ij}(a)=ae_{ij}-\bar ae_{ji}$, where $e_{ij}$ are the usual matrix units. These elements are subject to the following relations:
\begin{equation}\label{eq_relacionescolapsoSteinberg}
\begin{array}{l}
u_{ij}(a)=u_{ji}(-\bar a),\\
a\mapsto u_{ij}(a) \text{ is a linear map},\\
{[}u_{ij}(a),u_{jk}(b)]=u_{ik}(ab) \text{ for distinct $i$, $j$, $k$,}\\
{[}u_{ij}(a),u_{kl}(b)]=0\text{ for distinct $i$, $j$, $k$, $l$.}
\end{array}
\end{equation}
\end{Example}

Now, let $(A,-)$ be a unital algebra with involution, and let
$ \mathfrak{stu}_n(A,-)$ denote the Lie subalgebra generated by elements $u_{ij}(a)$, for $1\le i\ne j\le n$, $n\geq 3$, and $a\in A$, subject to the relations (\ref{eq_relacionescolapsoSteinberg}). Then the condition $u_{ij}(a)=0$ implies $a=0$ if and only if either $n\ge4$ and $A$ is associative or $n=3$ and $(A,-)$ is structurable. This Lie algebra is called the \emph{Steinberg unitary Lie algebra} by analogy with the Steinberg group in K-theory.

During the proof of the  previous result, Allison and Faulkner used  the following construction (also developed in \cite{paraSteinberg}).

Let $(A,-)$ be a structurable algebra. For $T\in\mathfrak{gl}(A)$, define $\bar T$ by $\bar T(x)=\overline{T(\bar x)}$. A set $T=(T_1,T_2,T_3)\in \mathfrak{gl}(A)^3$ is said to be a \emph{related triple} if
$$
\bar T_i(xy)=T_{i+1}(x)y+xT_{i+2}(y)
$$
for all $x,y\in A$ and for all $i=1,2,3$, where the subindices are taken modulo $3$. These triples form a Lie algebra denoted by $\mathfrak{trip}(A)$. A remarkable triple is the following:
\begin{equation}\label{eq_unrelatedtripleestandar}
\begin{array}{l}
T_i=L_{\bar b}L_{a}-L_{\bar a}L_{b},\\
T_{i+1}=R_{\bar b}R_{a}-R_{\bar a}R_{b},\\
T_{i+2}=R_{\bar ab-\bar ba}+L_{b}L_{\bar a}-L_{ a}L_{\bar b},
\end{array}
\end{equation}
for $a,b\in A$, where $L_a$ and $R_a$ denote the left and right multiplications by $a$ in $A$.

Consider now the vector space
$$
 \mathcal{U}(A,-):=\mathfrak{trip}(A)\oplus u_{12}(A)\oplus u_{23}(A)\oplus u_{31}(A)
$$
with anticommutative multiplication such that $\mathfrak{trip}(A)$ is a subalgebra and the following conditions hold:
$$
\begin{array}{l}
{[}T,u_{i,i+1}(a)]=u_{i,i+1}(T_{i+2}(a)),\\
{[}u_{i,i+1}(a),u_{i+1,i+2}(b)]=-u_{i+2,i}(\overline{ab}),\\
{[}u_{i,i+1}(a),u_{i,i+1}(b)]=T  \text{ as in (\ref{eq_unrelatedtripleestandar}),}
\end{array}
$$
for any $a,b\in A$ and $T=(T_1,T_2,T_3)\in\mathfrak{trip}(A)$. Then $\mathcal{U}(A,-)$ is a Lie algebra 
isomorphic to $\mathfrak{stu}_3(A,-)/\mathfrak{z}$, where $\mathfrak{z}$ is the center of $\mathfrak{stu}_3(A,-)$ (\cite[Theorem~4.3]{paraSteinberg}).
This algebra $\mathcal{U}(A,-)$ is simple 
 if and only if $(A,-)$ is so. 
We call this algebra the \emph{Steinberg construction} attached to the structurable algebra $A$.
Moreover, it turns out that $\mathcal{U}(A,-)$ is isomorphic to the Lie algebra given by the Kantor's construction $\mathfrak{Kan}(A,-)$ considered in the previous subsection.


\section{Description of the gradings}

\subsection{Gradings on   Hurwitz algebras}\label{subsec_gradingsonhurwitz}

Given a Hurwitz algebra $C$ with norm $q$, we can construct, for each $0\ne \alpha\in \FF$, a new unital algebra with involution, denoted by ${\CD}(C,\alpha)$, by means of the so called Cayley-Dickson doubling process. This is the algebra defined on $C\times C$ with the multiplication given by
$$
(a,b)(c,d)=(ac+\alpha\bar d b,da+b\bar c)
$$
and the quadratic form
$$
q(a,b)=q(a)-\alpha q(b).
$$
The resulting algebra ${\CD}(C,\alpha) $ is a Hurwitz algebra if and only if $C$ is associative. Thus, according to our list of Hurwitz algebras, ${\CD}(\FF,\alpha)$ is isomorphic to $\K$, ${\CD}(\K,\alpha)$ is isomorphic to $\Q$ and
${\CD}(\Q,\alpha)$ is isomorphic to $\C$. The algebra ${\CD}(C,\alpha) $ is always $\mathbb{Z}_2$-graded, with even part $\{(a,0)\mid a\in C\}$ and odd part $\{(0,b)\mid b\in C\}$. In particular $\K$ is $\mathbb{Z}_2$-graded,
$\Q$ is $\mathbb{Z}_2^2$-graded and $\C$ is $\mathbb{Z}_2^3$-graded.

On the other hand the standard basis of $\C$ is associated to a fine $\mathbb{ Z}^2$-grading, called the \emph{Cartan grading} on $\C$, which is given by
\begin{equation}\label{eq_graddeCdeCartan}
\begin{aligned}
\C_{(0,0)}  &=  \mathbb{F}e_1 \oplus \mathbb{F}e_2,&&\\
\C_{(1,0)}  &=   \mathbb{F}u_1, \quad&  \C_{(-1,0)} &= \mathbb{F}v_1,\\
\C_{(0,1)}  &=   \mathbb{F}u_2, &  \C_{(0,-1)} &= \mathbb{F}v_2,\\
\C_{(-1,-1)}  &=   \mathbb{F}u_3, \quad&  \C_{(1,1)} &= \mathbb{F}v_3.
\end{aligned}
\end{equation}
The subalgebra $\span{e_1,e_2,u_1,v_1}$ can be identified with $\Q$, which turns out to be $\mathbb{Z}$-graded, for
\begin{equation}\label{eq_graddeQdeCartan}
 \Q_{0}  =  \mathbb{F}e_1 \oplus \mathbb{F}e_2,\qquad  \Q_{1}  =   \mathbb{F}u_1, \qquad  \Q_{-1 } = \mathbb{F}v_1.
\end{equation}


\subsection{Gradings on symmetric composition algebras}\label{subsec_gradensymmcompalg}

Gradings on symmetric composition algebras were classified in \cite[Theorem~4.5]{gradingssymcomalg}.
Every group grading on  a Hurwitz algebra $C$  is a grading on $pC$ (since the involution preserves the homogeneous components), and the gradings on both the Hurwitz algebra $C$ and its para-Hurwitz counterpart coincide when  $C$ has dimension at least $4$.
In particular we have a $\mathbb{Z}_2$-grading on $\pK$, a $\mathbb{Z}_2^2$-grading and  a $\mathbb{Z}$-grading
on $\pQ$ and a $\mathbb{Z}_2^3$-grading and a   $\mathbb{Z}^2$-grading on $\pC$.

There is a remarkable $\mathbb{Z}_3$-grading in the case of dimension $2$ which does not come from a grading on the corresponding Hurwitz algebra, namely,
$$
\pK_{\bar{0}}=0,\qquad \pK_{\bar{1}}= \FF e_1\qquad  \pK_{\bar{2}}= \FF e_2,
$$
where $e_1$ and $e_2$ are the orthogonal idempotents $(1,0)$ and $(0,1)$ in $\K=\FF\times \FF$, which in $\pK$ satisfy
$e_1*e_1=e_2$ and $e_2*e_2=e_1$.

Also, a natural  $\mathbb{Z}_3^2$-grading appears on the pseudo-octonion algebra $\Ok=(\mathfrak{sl}_3(\FF),*,q)$,    determined by
$$
\Ok_{(\bar1,\bar0)}= \FF \left( \begin{array}{ccc}
1 & 0 & 0\\ 0 & \omega & 0\\ 0 & 0 & \omega^2
\end{array} \right)                   \quad        \text{and}\quad
\Ok_{(\bar0,\bar1)}= \FF \left( \begin{array}{ccc}
0 & 1 & 0\\ 0 & 0 & 1\\ 1 & 0 & 0
\end{array} \right)
$$
(recall that $\omega$ is a primitive cubic root of $1$).

\subsection{Gradings on Lie algebras induced from gradings on symmetric composition algebras}

If $S$ and $S'$ are two symmetric composition algebras, the Lie algebra $\mathfrak{g}(S,S')$ is always $\mathbb{Z}_2^2$-graded, for
\begin{equation}\label{eq_laZ22dadaporlaconstrucciong}
\begin{array}{ll}
\mathfrak{g}(S,S')_{(0,0)}= \mathfrak{tri}(S)\oplus\mathfrak{tri}(S'),\quad &\mathfrak{g}(S,S')_{(1,1)}= \iota_0(S\otimes S'),\\
\mathfrak{g}(S,S')_{(1,0)}=\iota_1(S\otimes S') ,&\mathfrak{g}(S,S')_{(0,1)}=\iota_2(S\otimes S').
\end{array}
\end{equation}
Moreover, if $S=\oplus_{g\in G} S_g$ is $G$-graded, and $S'=\oplus_{g\in G'} S'_g$ is $G'$-graded, these gradings can be combined with the one in Equation~(\ref{eq_laZ22dadaporlaconstrucciong}), thus obtaining a grading on
$\mathfrak{g}(S,S')$ by the group $\mathbb{Z}_2^2\times G\times G'$.

As a consequence, if we consider the gradings by the groups $\mathbb{Z}_2^r$ on the symmetric composition algebras as in Subsection~\ref{subsec_gradensymmcompalg}, we get gradings on:
\begin{description}
\item[\textbf{(6g1)}] $\mathfrak{e}_6=\mathfrak{g}(\pK,\pC)$ by the group $\mathbb{Z}_2^2\times \mathbb{Z}_2\times \mathbb{Z}_2^3=\mathbb{Z}_2^6$, of type $( 48,1,0,7  )$;
\item[\textbf{(7g1)}] $\mathfrak{e}_7=\mathfrak{g}(\pQ,\pC)$ by the group $\mathbb{Z}_2^2\times \mathbb{Z}_2^2\times \mathbb{Z}_2^3=\mathbb{Z}_2^7$, of type $(96,0,3,7   )$;
    \item[\textbf{(8g1)}] $\mathfrak{e}_8=\mathfrak{g}(\pC,\pC)$ by the group $\mathbb{Z}_2^2\times \mathbb{Z}_2^3\times \mathbb{Z}_2^3=\mathbb{Z}_2^8$, of type $(192,0,0,14   )$.
\end{description}

Moreover, there is a distinguished $\mathbb{Z}_3$-grading on $\mathfrak{g}(S,S')$, obtained as the eigenspace decomposition of an order $3$ automorphism $\Theta$ induced by the triality automorphisms $\vartheta$ and $\vartheta'$ of $S$ and $S'$ respectively. It is  given by
$$
\Theta\vert_{\mathfrak{tri}(S)}=\vartheta,\quad
\Theta\vert_{\mathfrak{tri}(S')}=\vartheta',\quad
\Theta(\iota_i(x\otimes x'))=\iota_{i+1}(x\otimes x').
$$
This $\mathbb{Z}_3$-grading on $\mathfrak{g}(S,S')$ induced by $\Theta$ can be combined with any $G$ and $G'$-gradings on $S$ and $S'$ respectively, to get
   a grading by the group $\mathbb{Z}_3\times G\times G' $ on $\mathfrak{g}(S,S')$. If this procedure is applied to the $\mathbb{Z}_3$-grading on $\pK$ and to the $\mathbb{Z}_3^2$-grading on $\Ok$, we get gradings on:
\begin{description}
\item[\textbf{(6g2)}] $\mathfrak{e}_6=\mathfrak{g}(\pK,\Ok)$ by the group $\mathbb{Z}_3\times \mathbb{Z}_3\times \mathbb{Z}_3^2=\mathbb{Z}_3^4$, of type $(72,0,2   )$;
    \item[\textbf{(8g2)}] $\mathfrak{e}_8=\mathfrak{g}(\Ok,\Ok)$ by the group $\mathbb{Z}_3\times \mathbb{Z}_3^2\times \mathbb{Z}_3^2=\mathbb{Z}_3^5$, of type $( 240,0,0,2  )$.
\end{description}
\noindent
 And if we combine some of the gradings above (related to the primes $2$ and $3$), we obtain gradings on:
\begin{description}
 \item[\textbf{(6g3)}] $\mathfrak{e}_6=\mathfrak{g}(\pK,\pC)$ by the group $\mathbb{Z}_3\times \mathbb{Z}_3\times \mathbb{Z}_2^3=\mathbb{Z}_2^3\times \mathbb{Z}_3^2$, of type $( 64,7  )$;
  \item[\textbf{(6g4)}] $\mathfrak{e}_6=\mathfrak{g}(\pK,\Ok)$ by the group $\mathbb{Z}_3\times \mathbb{Z}_2\times \mathbb{Z}_3^2=\mathbb{Z}_2\times \mathbb{Z}_3^3$, of type $(  26,26 )$;
\item[\textbf{(7g2)}] $\mathfrak{e}_7=\mathfrak{g}(\pQ,\Ok)$ by the group $\mathbb{Z}_3\times \mathbb{Z}_2^2\times \mathbb{Z}_3^2=\mathbb{Z}_2^2\times \mathbb{Z}_3^3$, of type $( 81,26  )$;
    \item[\textbf{(8g3)}] $\mathfrak{e}_8=\mathfrak{g}(\pC,\Ok)$ by the group $\mathbb{Z}_3\times \mathbb{Z}_2^3\times \mathbb{Z}_3^2=\mathbb{Z}_2^3\times \mathbb{Z}_3^3=\mathbb{Z}_6^3$, of type $( 182,33  )$.
\end{description}

Finally, we can also combine the above $\mathbb{Z}_3$-grading on $\mathfrak{g}(S,S')$ induced by $\Theta$ with  gradings on $S$ and $S'$ by infinite groups, namely, the $\mathbb{Z}$-grading on $\pQ$ and the $\mathbb{Z}^2$-grading on $\pC$ as in Equations~(\ref{eq_graddeQdeCartan}) and (\ref{eq_graddeCdeCartan}) respectively. Thus we get  gradings on:
\begin{description}
    \item[\textbf{(6g5)}] $\mathfrak{e}_6=\mathfrak{g}(\pK,\pC)$ by the group $\mathbb{Z}_3\times \mathbb{Z}_3\times \mathbb{Z}^2=\mathbb{Z}^2\times \mathbb{Z}_3^2$, of type $( 60,9  )$;
    \item[\textbf{(7g3)}] $\mathfrak{e}_7=\mathfrak{g}(\pQ,\Ok)$ by the group $\mathbb{Z}_3\times \mathbb{Z} \times \mathbb{Z}_3^2=\mathbb{Z} \times \mathbb{Z}_3^3$, of type $( 55,0,26  )$;
    \item[\textbf{(8g4)}] $\mathfrak{e}_8=\mathfrak{g}(\pC,\Ok)$ by the group $\mathbb{Z}_3\times \mathbb{Z}^2\times \mathbb{Z}_3^2=\mathbb{Z}^2\times \mathbb{Z}_3^3$, of type $( 168,1,26  )$.
\end{description}

Now, consider the  $\mathbb{Z}$-grading on $\mathfrak{g}(S,S')$  described in \cite[\S4.2]{Weylpreprint}
as follows:
Let $1$ and $1'$ be the unity elements of two Hurwitz algebras $C$ and $C'$. They become the  {paraunits} of the corresponding para-Hurwitz algebras $S=pC$ and $S'=pC'$.
Consider the inner derivation $d:=\ad(\iota_0(1\otimes 1'))$ of $\mathfrak{g}(S,S')$,  which is a semisimple derivation with eigenvalues
$\pm2,\pm1,0$. Thus, the eigenspace decomposition gives the following
$\mathbb{ Z}$-grading ($5$-grading) on $ \mathfrak{g}(S,S')$:
\begin{equation}\label{eq_laZgrad}
 \begin{aligned}
 \mathfrak{g}(S,S')_{\pm2}&=\Sigma_{\pm}(S_0,S'_0), \\
  \mathfrak{g}(S,S')_{\pm1}&=\nu_{\pm}(S\otimes S'), \\
   \mathfrak{g}(S,S')_0 &= t_{S_0, S_0}\oplus t_{S'_0,S'_0}\oplus \iota_0 (S_0\otimes S'_0) \oplus  \FF\iota_0 (1\otimes 1'),
 \end{aligned}
\end{equation}
where   $S_0$ and $S_0'$ denote the subspaces of zero trace elements in $C$ and $C'$ (here $S=pC$ and $S'=pC'$),
and where
$$
\begin{aligned}
\Sigma_{\pm}(y,y')&:=t_{1,y}+t'_{1,y'}\pm \ii \iota_0(y\otimes 1 + 1\otimes y'), \\
\nu_{\pm}(y\otimes y')&:=\iota_1(y\otimes y')\mp \ii \iota_2(\bar{y}\otimes \bar{y}'),
\end{aligned}
$$
 for all $y\in S$, $y'\in S'$, and for a fixed scalar $\ii \in  \FF$ such that  $\ii^2=-1$.
It is clear that this $\mathbb{Z}$-grading can be refined with gradings coming  from $S$ or $S'$.
In particular, when the $\mathbb{Z}_2^r$-gradings on $pC$ and $pC'$ are used, we get gradings on:
\begin{description}
    \item[\textbf{(6g6)}] $\mathfrak{e}_6=\mathfrak{g}(\pK,\pC)$ by the group $\mathbb{Z}\times \mathbb{Z}_2\times \mathbb{Z}_2^3=\mathbb{Z}\times \mathbb{Z}_2^4$, of type $(  57,0,7 )$;
    \item[\textbf{(7g4)}] $\mathfrak{e}_7=\mathfrak{g}(\pQ,\pC)$ by the group $\mathbb{Z}\times \mathbb{Z}_2^2\times \mathbb{Z}_2^3=\mathbb{Z}\times \mathbb{Z}_2^5$, of type $( 106,3,7  )$;
    \item[\textbf{(8g5)}] $\mathfrak{e}_8=\mathfrak{g}(\pC,\pC)$ by the group $\mathbb{Z}\times \mathbb{Z}_2^3\times \mathbb{Z}_2^3=\mathbb{Z}\times \mathbb{Z}_2^6$, of type $( 206,0,14  )$.
\end{description}


Next we recall the  $ \mathbb{Z}^4$-grading on $ \mathfrak{ g}(\pC,S')$ described in \cite[\S4.3]{Weylpreprint}. Take the canonical generators $a_1=(1,0,0,0)$, $a_2=(0,1,0,0)$, $g_1=(0,0,1,0)$, $g_2=(0,0,0,1)$  of the group $\mathbb{ Z}^4$, and write $a_0=-a_1-a_2, g_0=-g_1-g_2$. Set the degrees of the $ \mathbb{Z}^4$-grading as follows:
\begin{equation}\label{eq_laZ4gradeng}
\begin{array}{rcccl}
\deg\, \iota_i(e_1\otimes s) & = & a_i & = & -\deg\, \iota_i(e_2 \otimes s), \\
\deg\, \iota_i(u_i\otimes s) & = & g_i & = & -\deg\, \iota_i(v_i \otimes s), \\
\deg\, \iota_i(u_{i+1}\otimes s) & = & a_{i+2}+g_{i+1} & = & -\deg\, \iota_i(v_{i+1} \otimes s), \\
\deg\, \iota_i(u_{i+2}\otimes s) & = & - a_{i+1}+g_{i+2} & = & -\deg\, \iota_i(v_{i+2} \otimes s),
\end{array}
\end{equation}
where $s\in S'$ and $B=\{e_1, e_2, u_0, u_1, u_2, v_0, v_1, v_2\}$ is the standard basis of the algebra $\C$ described in Subsection~\ref{subsec_Compositionalgebras}. Also set $\deg(t)=(0,0,0,0)$ for all $t\in \mathfrak{tri}(S')$, and $\deg  t_{x,y}=\deg  \iota_0(x\otimes s) + \deg \iota_0(y\otimes s)$, if $x,y\in B$.
A straightforward computation shows that the $\mathbb{Z}^4$-grading on $\mathfrak{g}(\pC,S')$ provided by this assignment is compatible with any grading on $\pC$ and on $S'$. In particular, the following gradings are obtained:
\begin{description}
    \item[\textbf{(6g7)}] A grading on $\mathfrak{e}_6=\mathfrak{g}( \pC,\pK)$ by the group $\mathbb{Z}^4\times \mathbb{Z}_2$, of type $( 72,1,0,1  )$;
    \item[\textbf{(7g5)}] A grading on $\mathfrak{e}_7=\mathfrak{g}( \pC,\pQ)$ by the group $\mathbb{Z}^4\times \mathbb{Z}_2^2 $, of type $( 120,0,3,1  )$;
    \item[\textbf{(8g6)}] A grading on $\mathfrak{e}_8=\mathfrak{g}(\pC,\pC)$ by the group $\mathbb{Z}^4\times \mathbb{Z}_2^3 $, of type $(216,0,0,8)$.
\end{description}


Observe that a  $ \mathbb{Z}^3$-grading can be defined on $ \mathfrak{ g}(\pQ,S')$ and  a  $ \mathbb{Z}^2$-grading on $ \mathfrak{ g}(\pK,S')$, both of them inherited directly from the $ \mathbb{Z}^4$-grading   on $ \mathfrak{ g}(\pC,S')$ given by Equation~(\ref{eq_laZ4gradeng}). (Here  $\{e_1,e_2,u_1,v_1\}$ is a basis of $\pQ$ and $\{e_1,e_2\}$ is a basis of $\pK$.) The grading on $ \mathfrak{ g}(\pQ,S')$ is given by the following assignment of degrees:
\begin{equation*}\label{eq_laZ3gradeng}
\begin{array}{rcccl}
\deg\, \iota_i(e_1\otimes s) & = & a_i & = & -\deg\, \iota_i(e_2 \otimes s), \\
\deg\, \iota_1(u_1\otimes s) & = & (0,0,1) & = & -\deg\, \iota_1(v_1 \otimes s), \\
\deg\, \iota_2(u_{1}\otimes s) & = & (0,1,1) & = & -\deg\, \iota_2(v_{1} \otimes s), \\
\deg\, \iota_0(u_{1}\otimes s) & = & (1,1,1) & = & -\deg\, \iota_0(v_{1} \otimes s),\\
\deg(t_{e_1,e_2})&=&(0,0,0)&=&\deg(t_{u_1,v_1}),\\
\deg(t_{e_1,u_1})&=&(-1,0,1)&=&-\deg(t_{e_2,v_1}),\\
\deg(t_{e_1,v_1})&=&(-1,-2,-1)&=&-\deg(t_{e_2,u_1}),
\end{array}
\end{equation*}
 for  $a_1=(1,0,0)$, $a_2=(0,1,0)$, $a_0=(-1,-1,0)$, with $\deg(\mathfrak{tri}(S'))=(0,0,0)$.
The grading on $ \mathfrak{ g}(\pK,S')$ is given by:
\begin{equation*}
\begin{array}{rcccl}
\deg\, \iota_1(e_1\otimes s) & = & (1,0) & = & -\deg\, \iota_1(e_2 \otimes s), \\
\deg\, \iota_2(e_{1}\otimes s) & = & (0,1) & = & -\deg\, \iota_2(e_{2} \otimes s), \\
\deg\, \iota_0(e_{1}\otimes s) & = & (-1,-1) & = & -\deg\, \iota_0(e_{2} \otimes s),
\end{array}
\end{equation*}
with $\deg(\mathfrak{tri}(\pK))=(0,0)= \deg(\mathfrak{tri}(S'))$.  Again these gradings can be combined with the $\mathbb{Z}_2^3$-grading on the symmetric composition algebra $\pC$ to get:
\begin{description}
    \item[\textbf{(6g8)}] A grading on $\mathfrak{e}_6=\mathfrak{g}( \pK,\pC)$ by the group $\mathbb{Z}^2\times \mathbb{Z}_2^3$, of type $(48,1,0,7   )$; 
    \item[\textbf{(7g6)}] A grading on $\mathfrak{e}_7=\mathfrak{g}( \pQ,\pC)$ by the group $\mathbb{Z}^3\times \mathbb{Z}_2^3 $, of type $( 102,0,1,7  )$.
\end{description}

\subsection{Gradings on some Jordan  algebras}\label{Subsec_gradingsonJordanofdegree4}

We are going to describe some gradings on the Jordan algebra  $H_4(C,*)=\{x=(x_{ij})\in \Mat_{4\times 4}(C)\mid x_{ij}=\bar x_{ji}\}$ for $C$ some associative Hurwitz algebra (that is, up to isomorphism, $C\in\{\FF,\K,\Q\}$). The reader may consult \cite[Chapter 5]{librogradings} for the description of gradings on simple Jordan algebras.

Observe first that the Kronecker product gives an isomorphism of associative algebras:
$$
\Mat_{2\times 2}(\FF)\otimes \Mat_{2\times 2}(\FF)\to \Mat_{4\times 4}(\FF),\quad
a\otimes b\mapsto a\otimes b=\left(\begin{array}{cc}a_{11}b&a_{12}b\\a_{21}b&a_{22}b\end{array}\right).
$$
Also, there is the natural isomorphism of associative algebras,
\begin{equation}\label{eq_identificacionMat4Q}
\Mat_{4\times 4}(\FF)\otimes C\cong\Mat_{4\times 4}(C),\quad (a_{ij})\otimes x\mapsto (a_{ij}x).
\end{equation}
As $\Mat_{2\times 2}(\FF)$ is isomorphic to $\Q$, it inherits a $\mathbb{Z}$-grading and a $\mathbb{Z}_2^2$-grading, so that the previous identifications allow us to define gradings on
$\Mat_{4\times 4}(\FF)$ by the groups $\mathbb{Z}_2^2\times \mathbb{Z}_2^2$ and $\mathbb{Z}\times\mathbb{Z}_2^2$,
on $\Mat_{4\times 4}(\K)$ by the groups $\mathbb{Z}_2^2\times \mathbb{Z}_2^2\times \mathbb{Z}_2 $ and $\mathbb{Z}\times\mathbb{Z}_2^2\times\mathbb{Z}_2$,
and on $\Mat_{4\times 4}(\Q)$ by the groups $\mathbb{Z}_2^2\times \mathbb{Z}_2^2\times \mathbb{Z}_2^2 $ and $\mathbb{Z}\times\mathbb{Z}_2^2\times\mathbb{Z}_2^2$.
Trivially any grading on  the associative algebra $\Mat_{4\times 4}(C)$ is a grading of the Jordan algebra $\Mat_{4\times 4}(C)^+$.
The point is that, for the previously described gradings, the Jordan subalgebra $H_4(C,*)\le \Mat_{4\times 4}(C)^+$ is a graded subspace, so that:
\begin{itemize}
    \item[$\bullet$ ] $H_4(\FF,*)$ is $\mathbb{Z}_2^4$ and $\mathbb{Z}\times\mathbb{Z}_2^2$-graded;
    \item[$\bullet$ ] $H_4(\K,*)$ is $\mathbb{Z}_2^5$ and $\mathbb{Z}\times\mathbb{Z}_2^3$-graded; 
    \item[$\bullet$ ]  $H_4(\Q,*)$ is $\mathbb{Z}_2^6$ and $\mathbb{Z}\times\mathbb{Z}_2^4$-graded.
\end{itemize}
Let us explain this with some extra detail. Let us denote by $ q_0$ the identity matrix of degree $2$, and consider the matrices
\begin{equation}\label{eq_laZ22decuaternios}
q_1=\left(\begin{array}{cc}0&1\\1&0\end{array}\right),\quad
q_2=\left(\begin{array}{cc}1&0\\0&-1\end{array}\right),\quad
q_3=\left(\begin{array}{cc}0&-1\\1&0\end{array}\right)=q_1q_2.
\end{equation}
Then the assignment $\deg(q_1)=(\bar1,\bar0)$ and  $\deg(q_2)=(\bar0,\bar1)$ gives the $\mathbb{Z}_2^2$-grading on $\Mat_{2\times 2}(\FF)$. The $\mathbb{Z}_2^4$-grading on $\Mat_{4\times 4}(\FF)$
has $16$ one-dimensional homogeneous components, where $q_i\otimes q_j$ has degree $(\deg( q_i),\deg( q_j))$. The subset of homogeneous elements
  $\{ q_i\otimes  q_j, q_3\otimes q_3\mid i,j=0,1,2\}$ spans $H_4(\FF,*)$, and hence the $10$-dimensional space $H_4(\FF,*)$ is also $\mathbb{Z}_2^4$-graded. Moreover, as the $\mathbb{Z}_2$-grading on $\K$ is given by $\K_{\bar0}=\span{1=e_1+e_2}$ and $\K_{\bar1}=\span{e_1-e_2}$, then $H_4(\K,*)\subset \Mat_{2\times 2}(\FF)\otimes\Mat_{2\times 2}(\FF)\otimes\K$ is spanned by the following subset of homogeneous elements for the $\mathbb{Z}_2^5$-grading:
  $\{ q_i\otimes  q_j\otimes 1, q_3\otimes q_3\otimes1\mid i,j=0,1,2\}
  \cup
  \{ q_i\otimes  q_3\otimes (e_1-e_2), q_3\otimes q_i\otimes(e_1-e_2)\mid i=0,1,2\}$. The remaining cases are dealt with in the same way.

\subsection{Gradings on Lie algebras obtained from Kantor's and Steinberg's constructions}

Recall that we can get the exceptional Lie algebras of the $E$ series by means of  Kantor's construction applied to  the  structurable algebras  $ {\CD}(H_4(C,*))$, for an associative Hurwitz algebra $C$. In turn, these structurable algebras are obtained from the Jordan algebras $J=H_4(C,*)$ by means of the Cayley-Dickson doubling process explained in Example~\ref{ex_CDprocesoparaestructurables}. This doubling process provides a $\mathbb{Z}_2$-grading
as usual, with even part $J$ and odd part $vJ$, clearly compatible with any grading on $J$. At the same time, any $G$-grading on a structurable algebra $(A,-)$ provides a $\mathbb{Z}\times G$-grading on $\mathfrak{Kan}(A,-)$ and a
$\mathbb{Z}_2^2\times G$-grading on $\mathcal{U}(A,-)$ (which is isomorphic to $\mathfrak{Kan}(A,-)$). Thus we have another source of gradings on our Lie algebras. If Kantor's construction is applied to $ {\CD}(H_4(C,*))$, and the $\mathbb{Z}_2$-grading induced by the Cayley-Dickson doubling process and the finite gradings on Subsection~\ref{Subsec_gradingsonJordanofdegree4} are combined, we obtain gradings on:
\begin{description}
    \item[\textbf{(6g9)}]   $\mathfrak{e}_6=\mathfrak{Kan}( {\CD}(H_4(\FF)))$ by the group $\mathbb{Z} \times \mathbb{Z}_2^5$, of type $( 73,0,0,0,1  )$;
    \item[\textbf{(7g7)}]   $\mathfrak{e}_7=\mathfrak{Kan}( {\CD}(H_4(\K)))$ by the group $\mathbb{Z} \times \mathbb{Z}_2^6 $, of type $( 127,0,0,0,0,1  )$;
    \item[\textbf{(8g7)}]    $\mathfrak{e}_8=\mathfrak{Kan}( {\CD}(H_4(\Q)))$ by the group $\mathbb{Z} \times \mathbb{Z}_2^7 $, of type $( 241,0,0,0,0,0,1  )$.
 \end{description}

In the same vein, but using the infinite gradings on       Subsection~\ref{Subsec_gradingsonJordanofdegree4}, we get gradings on:
\begin{description}
    \item[\textbf{(6g10)}]   $\mathfrak{e}_6=\mathfrak{Kan}( {\CD}(H_4(\FF)))$ by the group $\mathbb{Z}^2 \times \mathbb{Z}_2^3$, of type $( 60,7,0,1  )$;
    \item[\textbf{(7g8)}]   $\mathfrak{e}_7=\mathfrak{Kan}( {\CD}(H_4(\K)))$ by the group $\mathbb{Z}^2 \times \mathbb{Z}_2^4 $, of type $( 102,13,0,0,1  )$;   
     \item[\textbf{(8g8)}]   $\mathfrak{e}_8=\mathfrak{Kan}( {\CD}(H_4(\Q)))$ by the group $\mathbb{Z}^2 \times \mathbb{Z}_2^5 $, of type $( 180,31,0,0,0,1  )$.
\end{description}

Moreover, if   we use Steinberg's construction applied to $ {\CD}(H_4(C,*))$, the $\mathbb{Z}_2$-grading induced by the Cayley-Dickson doubling process and the finite gradings in Subsection~\ref{Subsec_gradingsonJordanofdegree4} can be combined to get gradings on:
\begin{description}
    \item[\textbf{(6g11)}]   $\mathfrak{e}_6=\mathcal{U}( {\CD}(H_4(\FF)))$ by the group $ \mathbb{Z}_2^7$, of type $( 72,0,0,0,0,1  )$;
    \item[\textbf{(7g9)}]   $\mathfrak{e}_7=\mathcal{U}( {\CD}(H_4(\K)))$ by the group $ \mathbb{Z}_2^8 $, of type $( 126,0,0,0,0,0,1  )$;
    \item[\textbf{(8g9)}]    $\mathfrak{e}_8=\mathcal{U}( {\CD}(H_4(\Q)))$ by the group $ \mathbb{Z}_2^9 $, of type $( 240,0,0,0,0,0,0,1  )$.
\end{description}
It is not difficult to see that the $\mathbb{Z} \times \mathbb{Z}_2^5$-grading (respectively $\mathbb{Z} \times \mathbb{Z}_2^6$ and $\mathbb{Z} \times \mathbb{Z}_2^7$) obtained in  $\mathfrak{e}_6=\mathcal{U}( {\CD}(H_4(\FF)))$
 (respectively $\mathfrak{e}_7=\mathcal{U}( {\CD}(H_4(\K)))$ and $\mathfrak{e}_8=\mathcal{U}( {\CD}(H_4(\Q)))$) is isomorphic to the grading \textbf{(6g9)}
  (respectively \textbf{(7g7)} and \textbf{(8g7)}).
\smallskip

\begin{Remark}
As we know \cite{e6} about the existence of a $\mathbb{Z}_4 \times \mathbb{Z}_2^4$-grading on $\mathfrak{e}_6$, of type $( 48,13,0,1  )$, we would like to find a
  $\mathbb{Z}_4$-grading on the Lie algebra obtained by means of Kantor's construction attached to a  structurable algebra. That can be done for    $A={\CD}(H_4(C))=J\oplus vJ$, with $J=H_4(C,*)$   the Jordan algebra of hermitian matrices with coefficients in an associative Hurwitz algebra $C$. In such a case, $S=\FF v$, and it is an straightforward computation, taking into account Equations~(\ref{eq_productoenKantor}) and (\ref{eq_CDestructurable}), that $L=\mathfrak{Kan}(A)$ is $\mathbb{Z}_4$-graded as follows:
\begin{equation*}\label{eq_laZ4estructurableCD}
\begin{array}{l}
L_{\bar0}= \FF v{\tilde{\ }}\oplus (V_{J,J}+V_{vJ,vJ})\oplus   \FF v   ,   \\
L_{\bar1}= J\oplus (vJ){\tilde{\ }} ,  \\
L_{\bar2}=  V_{J,vJ},  \\
L_{\bar3}=J{\tilde{\ }}\oplus vJ.
\end{array}
\end{equation*}
This grading is compatible with  any grading  on $J$, so that we can combine it with the gradings described in Subsection~\ref{Subsec_gradingsonJordanofdegree4} to get
gradings on
   $\mathfrak{e}_6=\mathfrak{Kan}( {\CD}(H_4(\FF)))$ by the group $\mathbb{Z}_4 \times \mathbb{Z}_2^4$, on   $\mathfrak{e}_7=\mathfrak{Kan}( {\CD}(H_4(\K)))$ by the group $\mathbb{Z}_4 \times \mathbb{Z}_2^5 $, and on   $\mathfrak{e}_8=\mathfrak{Kan}( {\CD}(H_4(\Q)))$ by the group $\mathbb{Z}_4 \times \mathbb{Z}_2^6 $.
   But these are not fine!, so we must continue the search in order to explain several gradings by groups with factors $\mathbb{Z}_4$.
   \end{Remark}

\subsection{A $\mathbb{Z}_4\times\mathbb{Z}_2^3$-grading on the Jordan   algebra $H_4(\Q)$ and related gradings on the exceptional Lie algebras}\label{subsec_laotragradcon4quenecesitamos}

A  graded division (associative) algebra $D$ is a graded algebra such that every homogeneous element is invertible. If the support of such a grading is $H$ and $G$ is a group
  containing $H$ as a subgroup, and we have a $G$-graded right $D$-module
$V$ (that is, $V_gD_h \subset V_{g+h}$ for any $g\in G$ and $h\in H$), then the
division property of $D$ forces $V$ to be a free right $D$-module containing bases consisting
of homogeneous elements, according to \cite[\S2]{Albertoclasicas}. Then we have
a $G$-grading induced on $R = \End_D(V )$ given by
$f\in R_g$ if $f(V_{g'}) \subset V_{g+g'}$
for any $g'\in G$.

Let $\tau\colon D\to D$ be a graded antiautomorphism, that is, $\tau$ is an antiautomorphism
with $\tau(D_h)=D_h$ for any $h\in H$ (which implies that necessarily $\tau$ is an involution, that is $\tau^2=\id_D$). Let $b\colon V\times V\to D$ a sesquilinear form
($b$ is $\FF$-bilinear, $b(v_1,v_2)=\tau(b(v_2,v_1))$ and $b(v_1,v_2d)=b(v_1,v_2)d$ for any $v_1,v_2\in V$ and $d\in D$)
compatible with the grading, that is, $b(V_g,V_{g'})\subset D_{g+g'}$. Let $*$ be the adjoint relative to this form   ($b(f(v_1),v_2)=b(v_1,f^*(v_2))$ if $f\in R$). The point is that the sets of hermitian and skew-hermitian elements 
  $H(R,*)=\{f\in R\mid f^*=f\}$ and $K(R,*)=\{f\in R\mid f^*=-f\}$ are graded subspaces. Moreover, it is proven in \cite{Albertoclasicas} that essentially all the gradings in $K(R,*)$ and in  $H(R,*)$  are obtained in this way.

Note that the quaternion algebra $\Q= \Mat_{2\times2}(\FF) $ is a $\mathbb{Z}_2^2$-graded division algebra with the grading given by the matrices in Equation~(\ref{eq_laZ22decuaternios}),
\begin{equation}\label{eq_laZ22decuaterniosESTAVEZSI}
\Q_{(\bar0,\bar0)}=\FF 1,\quad
\Q_{(\bar1,\bar0)}=\FF q_1,\quad
\Q_{(\bar0,\bar1)}=\FF q_2,\quad
\Q_{(\bar1,\bar1)}=\FF q_3.
\end{equation}
 There are two involutions compatible with this grading.
 The involution  $\tau^o$ given by
 $$
 q_1^{\tau^o}=q_1,\quad
 q_2^{\tau^o}=q_2,\quad
 q_3^{\tau^o}=-q_3,\quad
 $$
is the usual transpose involution (an orthogonal involution); while the involution $\tau^s=-$ given
by
$$
 \bar q_i =-q_i,\quad \forall i=1,2,3,
 $$
is the standard conjugation of the quaternion algebra $\Q$ (a symplectic
involution).

Take $(D=\Q,\tau^o)$ as above but with the following grading, equivalent to (\ref{eq_laZ22decuaterniosESTAVEZSI}),
\begin{equation*}
\Q_{(\bar0,\bar0)}=\FF 1,\quad
\Q_{(\bar2,\bar0)}=\FF q_1,\quad
\Q_{(\bar0,\bar1)}=\FF q_2,\quad
\Q_{(\bar2,\bar1)}=\FF q_3.
\end{equation*}
Take $B=\{v_0,v_1\}$ a homogeneous $D$-basis in $V$, a graded right free $D$-module of dimension 2,
with $\deg(v_0)=(\bar0,\bar0)$ and  $\deg(v_1)=(\bar1,\bar0)$. We have chosen the degrees such that the
   sesquilinear form   $b\colon V\times V\to D$ given by the matrix
   $A=\left(\begin{array}{cc}1&0\\0&q_1\end{array}\right) $ relative to $B$ is
   compatible with the grading (since $2\deg(v_0)=\deg(1)$ and $2\deg(v_1)=\deg(q_1)$).
   Now for $x=\left(\begin{array}{cc}p_1&p_2\\p_3&p_4\end{array}\right) \in\Mat_{2\times2}(\Q)\simeq\End_{\Q}(V)=R$,
   we have $x^*=A^{-1}(\tau^o(x))^tA=
   \left(\begin{array}{cc}\tau^o(p_1)&\tau^o(p_3)q_1\\-q_1\tau^o(p_2)&-q_1\tau^o(p_4)q_1\end{array}\right)$,
   and hence
   $$
   K(\Mat_{2\times2}(\Q),*)=\left\{\left(\begin{array}{cc}\alpha q_3&-\tau^o(p)q_1\\p&\beta q_2\end{array}\right) \mid  \alpha,\beta\in\FF,p\in\Q \right\}=: K
   $$
inherits the $\mathbb{Z}_4\times\mathbb{Z}_2$-grading,   with $6$ pieces of dimension one,
$$
K=K_{(\bar2,\bar1)}\oplus K_{(\bar0,\bar1)}\oplus K_{(\bar1,\bar0)}\oplus K_{(\bar3,\bar0)}\oplus K_{(\bar1,\bar1)}\oplus K_{(\bar3,\bar1)},
$$
and also, for $H:=H(\Mat_{2\times2}(\Q),*)$,
$$
   H=\left\{\left(\begin{array}{cc}p_1&\tau^o(p)q_1\\p&p_2\end{array}\right) \mid  p_1\in\span{1,q_1,q_2},p_2\in\span{1,q_1,q_3},p\in\Q \right\}
   $$
inherits the $\mathbb{Z}_4\times\mathbb{Z}_2$-grading,
$$
H=H_{(\bar2,\bar1)}\oplus H_{(\bar0,\bar1)}\oplus H_{(\bar1,\bar0)}\oplus H_{(\bar3,\bar0)}\oplus H_{(\bar1,\bar1)}\oplus H_{(\bar3,\bar1)}\oplus H_{(\bar0,\bar0)}\oplus H_{(\bar2,\bar0)},
$$
of type $(6,2)$ since all the above homogeneous components are one-dimensional except for two of them, namely,
$$
\small{
H_{(\bar0,\bar0)}= \span{\left(\begin{array}{cc}1&0\\0&0\end{array}\right) ,\left(\begin{array}{cc}0&0\\0&1\end{array}\right) },  \quad
 H_{(\bar2,\bar0)}=\span{\left(\begin{array}{cc}q_1&0\\0&0\end{array}\right) ,\left(\begin{array}{cc}0&0\\0&q_1\end{array}\right) } }.
$$

Next we identify $\Mat_{4\times4}(\Q)$ with $\Mat_{2\times2}(\FF)\otimes \Mat_{2\times2}(\FF)\otimes\Q$ and hence with
$\Q\otimes\Mat_{2\times2}(\Q)=\Q\otimes R$, and consider here the (symplectic) involution given by $\tau^s\otimes *$.
The Jordan algebra $J=H_4(\Q)$ lives here as $J=\{q\otimes x\in\Q\otimes R \mid \tau^s(q)\otimes x^*=q\otimes x\}$, which can be identified with $K(\Q,\tau^s)\otimes K(\Mat_{2\times2}(\Q),*)\oplus H(\Q,\tau^s)\otimes H(\Mat_{2\times2}(\Q),*)$. In this way, by combining the $\mathbb{Z}_2^2$-grading on $\Q$ given by Equation~(\ref{eq_laZ22decuaterniosESTAVEZSI}) with the above $\mathbb{Z}_4\times\mathbb{Z}_2$-grading on $R$, we get a $\mathbb{Z}_4\times\mathbb{Z}_2^3$-grading on $J=H_4(\Q)$ of type $3(6,0)+1(6,2)=(24,2)$. Also we get a $\mathbb{Z}_4\times\mathbb{Z}_2^3$-grading on $K(\Q\otimes R,\tau^s\otimes *)\cong\Der(H_4(\Q))$  (a Lie algebra of type $\mathfrak{c}_4$) of type $1(6,0)+3(6,2)=(24,6)$.\medskip

Of course this grading on $J$ induces a   $\mathbb{Z}_4\times\mathbb{Z}_2^4$-grading on the structurable algebra ${\CD}(H_4(\Q))$, of type $(48,4)$. Now note that, according to \cite{Allisonparaderystr},
\begin{equation}\label{eq_e678ylaestructurablededim56}
\Der({\CD}(H_4(\Q)),-)\simeq\mathfrak{e}_6,\quad
\mathfrak{Instr}({\CD}(H_4(\Q)))\simeq\mathfrak{e}_7,\quad
\mathcal{U}({\CD}(H_4(\Q)))\simeq\mathfrak{e}_8.
\end{equation}
In particular, every $G$-grading on $ {\CD}(H_4(\Q))$ induces a $G$-grading on $\mathfrak{e}_6$, a $G\times \mathbb{Z}_2$-grading on $\mathfrak{e}_7$ and a $G\times \mathbb{Z}_2^2$-grading on $\mathfrak{e}_8$.
In our case we get
gradings on:
\begin{description}
    \item[\textbf{(6g12)}]   $\mathfrak{e}_6=\Der( {\CD}(H_4(\Q)))$ by the group $\mathbb{Z}_4 \times \mathbb{Z}_2^4$, of type $( 48,13,0,1  )$;
    \item[\textbf{(7g10)}]   $\mathfrak{e}_7=\mathfrak{Instr}( {\CD}(H_4(\Q)))$ by the group $\mathbb{Z}_4 \times \mathbb{Z}_2^5 $, of type $( 98,15,0,0,1  )$;
    \item[\textbf{(8g10)}]    $\mathfrak{e}_8=\mathcal{U}( {\CD}(H_4(\Q)))$ by the group $\mathbb{Z}_4 \times \mathbb{Z}_2^6 $, of type $( 192,25,0,0,0,1 )$.
\end{description}


\subsection{A $\mathbb{Z}_4^3$-grading on the structurable  algebra ${\CD}(H_4(\Q))$}\label{subsec_Z43enlaestructurable}

There is a $\mathbb{Z}_4^3$-grading on the structurable algebra $\mathcal{A}={\CD}(H_4(\Q))$ which is not explained in terms of  the Cayley-Dickson process.
This is a very interesting grading in which every nonzero homogeneous component is one-dimensional. Let us describe it. The information is extracted from   \cite{Z43proceedings}.

Identify, as in Equation~(\ref{eq_identificacionMat4Q}),
 $\Mat_{4\times 4}(\Q)$ with
$\Mat_{4\times 4}(\FF)\otimes \Q$. The involution $(q_{ij})^*=(\bar q_{ji})$ in  $\Mat_{4\times 4}(\Q)$
is, under such correspondence,  the
tensor product of the matrix transpose on $\Mat_{4\times 4}(\FF)$ and the standard involution on $\Q$ (under the identification of $\Q$ with  $\Mat_{2\times 2}(\FF)$, this involution acts as follows: $\bar E_{11}=E_{22}$,
$\bar E_{22}=E_{11}$, $\bar E_{12}=-E_{12}$, $\bar E_{21}=-E_{21}$). In particular,
  the Jordan subalgebra of symmetric elements
$H_4(\Q) = \{q=(q_{ij})\in \Mat_{4\times 4}(\Q)\mid q^*=q\}$ is identified with
$
H_4(\FF)\otimes\span{E_{11}+E_{22}}\oplus
K_4(\FF)\otimes\span{E_{11}-E_{22},E_{12},E_{21}}$, where $H_4(\FF)$ and  $K_4(\FF)$ denote the subspaces of symmetric and skewsymmetric matrices of size 4, respectively,
and hence with
$$
\mathcal{J}=\left\{\left(\begin{array}{cc}z&x\\y&z^t\end{array}\right)\mid x=-x^t,y=-y^t,\,x,y,z\in \Mat_{4\times 4}(\FF)\right\}.
$$
Consider the following $\mathbb{Z}_4$-grading on the structurable algebra $\mathcal{A}=\mathcal{J}\oplus v\mathcal{J}$:
\begin{equation}\label{eq_laZ4delaCDH4Q}
\begin{array}{l}
\mathcal{A}_{\bar0}=\left\{\left(\begin{array}{cc}z&0\\0&z^t\end{array}\right)\mid z\in \Mat_{4\times 4}(\FF)\right\},\\
\mathcal{A}_{\bar1}=\left\{\left(\begin{array}{cc}0&x\\0&0\end{array}\right)+
v\left(\begin{array}{cc}0&0\\y&0\end{array}\right)\mid x,y \in K_4(\FF)\right\},\\[4pt]
\mathcal{A}_{\bar2}=v\mathcal{A}_{\bar0},\\[4pt]
\mathcal{A}_{\bar3}=\left\{v\left(\begin{array}{cc}0&x\\0&0\end{array}\right)+
\left(\begin{array}{cc}0&0\\y&0\end{array}\right)\mid x,y \in K_4(\FF)\right\}.
\end{array}
\end{equation}
Note that if $u\in  \textrm{GL}(4,\FF)$ is an invertible matrix, we can consider the following automorphism of $H_4(\Q)$,
$$
\Psi(u):\left(\begin{array}{cc}z&x\\y&z^t\end{array}\right)\mapsto \left(\begin{array}{cc}uzu^{-1}&uxu^t\\(u^{-1})^tyu^{-1}&(uzu^{-1})^t \end{array}\right),
$$
which extends to $\mathcal{A}$ in a natural way (also denoted $\Psi(u)$). For any natural number $n$, let $\xi$ be a primitive $n$th root of $1$ and consider the following matrices
{\small
\begin{equation}\label{eq_Paulimatrices}
P_n=\left(\begin{array}{ccccc}1&0&\dots&\dots&0\\0&\xi&0&\dots&0
\\\vdots&\ddots&\ddots&\ddots&\vdots\\0&\ldots&0&\xi^{n-2}&0\\0&\ldots&\ldots&0&\xi^{n-1}\end{array}\right),
\quad
Q_n=\left(\begin{array}{ccccc}0&1&0&\dots&0\\0&0&1&\dots&0
\\{\vdots}&\ddots&\ddots&\ddots&\vdots\\0&\ldots&0&0&1\\1&0&\ldots&\ldots&0\end{array}\right),
\end{equation}
}%
also called Pauli matrices ($P_n$ sometimes will be denoted by $P_\xi$).
Take now the Pauli matrices $X=P_{4}$ and $Y=Q_4$ and note that then $\Psi(X)$ and $\Psi(Y)$ are order four automorphisms of $\mathcal{A}$ (which neither commute nor  anticommute).

Also, consider for any skew-symmetric $4\times 4$ matrix, its Pfaffian adjoint $\hat x$:
$$
x=\left(\begin{array}{cccc}0&\alpha&\beta&\gamma\\-\alpha&0&\delta&\varepsilon\\-\beta&-\delta&0&\zeta\\
-\gamma&-\varepsilon&-\zeta&0\end{array}\right)\in K_4(\FF),\quad
\hat{x}=\left(\begin{array}{cccc}0&-\zeta&\varepsilon&-\delta\\ \zeta&0&-\gamma&\beta\\-\varepsilon&\gamma&0&-\alpha\\
\delta&-\beta&\alpha&0\end{array}\right).
$$
(Note that this differs from \cite{Z43proceedings}, where $-\hat x$ is considered.)
Now consider the order $4$ automorphism   $\pi\colon\mathcal{A}\to\mathcal{A}$, whose restriction to $\mathcal{A}_{\bar0}\oplus\mathcal{A}_{\bar2}$ is the identity, and such that:
$$
\pi\left(\begin{array}{cc}0&x\\y&0\end{array}\right)=v\left(\begin{array}{cc}0&-\hat y\\ \hat x&0\end{array}\right),\quad
\pi\left(v\left(\begin{array}{cc}0&x\\y&0\end{array}\right)\right)=\left(\begin{array}{cc}0&-\hat y\\ \hat x&0\end{array}\right),
$$
on $\mathcal{A}_{\bar1}\oplus\mathcal{A}_{\bar3}$.
If $\xi\in \FF$ is chosen such that $\xi^2=\ii$, then $\pi\Psi(\xi X)$ and $\Psi(Y)$ are two order $4$ commuting automorphisms that preserve the $\mathbb{Z}_4$-grading given by Equation~(\ref{eq_laZ4delaCDH4Q}), and a  $\mathbb{Z}_4^3$-grading on $\mathcal{A}$ is obtained whose homogeneous components are the intersection of the homogeneous components of the $\mathbb{Z}_4$-grading with the common eigenspaces for $\pi\Psi(\xi X)$ and $\Psi(Y)$.

\begin{Remark}
This structurable algebra of dimension $56$ (the only simple one of such dimension) has a model which is better known. It is defined on the vector space
$$
\left(\begin{array}{cc}
\FF&\mathbb{A}\\
\mathbb{A}&\FF\end{array}\right)
$$
with multiplication given by
$$
\left(\begin{array}{cc}
\alpha_1&x_1\\
x_1'&\beta_1\end{array}\right)
\left(\begin{array}{cc}
\alpha_2&x_2\\
x_2'&\beta_2\end{array}\right)=
\left(\begin{array}{cc}
\alpha_1\alpha_2+T(x_1,x_2')&\alpha_1x_2+\beta_2x_1+x_1'\times x_2'\\
\alpha_2x_1'+\beta_1x_2'+x_1\times x_2&\beta_1\beta_2+T(x_2,x_1')\end{array}\right),
$$
where $T$ denotes the map $T\colon \mathbb{A}\times \mathbb{A}\to\FF$ given by $T(x,y)=T(xy)$
and $\times$ denotes the so called Freudenthal cross product   defined by $T(x\times y,z)=N(x,y,z)$ if $x,y,z\in\mathbb{A}$.
The involution is given by
$$
\overline{\left(\begin{array}{cc}
\alpha &x \\
x '&\beta \end{array}\right)}=\left(\begin{array}{cc}
\beta &x \\
x '&\alpha \end{array}\right).
$$


Although this is isomorphic (as an algebra with involution), to ${\CD}(H_4(\Q))$, it was previously studied as an example of Brown algebra. Garibaldi \cite{Garibaldi} discusses the connections between this algebra and the groups of types $E_6$ and $E_7$.
\end{Remark}

\subsection{More gradings on the exceptional Lie algebras by groups with factors $\mathbb{Z}_4$}\label{subsec_gradingsconcuatros}

The $\mathbb{Z}_4^3$-grading above on ${\CD}(H_4(\Q))$ immediately induces the following gradings on:
\begin{description}
\item[\textbf{(8g11)}]    $\mathfrak{e}_8=\mathfrak{Kan}( {\CD}(H_4(\Q)))$ by the group $\mathbb{Z} \times \mathbb{Z}_4^3 $, of type $( 123,40,15 )$;
  \item[\textbf{(8g12)}]    $\mathfrak{e}_8=\mathcal{U}( {\CD}(H_4(\Q)))$ by the group $\mathbb{Z}_2^2 \times \mathbb{Z}_4^3 $, of type $(216,14,0,1)$;
\end{description}
if we take into consideration the $\mathbb{Z}$-grading provided by Kantor's construction and the $\mathbb{Z}_2^2$-grading provided by   Steinberg's construction.\smallskip

Furthermore, recall that the algebras $\mathfrak{e}_6$ and $\mathfrak{e}_7$ can be obtained too from the structurable algebra $ {\CD}(H_4(\Q))$ as in (\ref{eq_e678ylaestructurablededim56}).   Hence we get also gradings on:
\begin{description}
\item[\textbf{(6g13)}]    $\mathfrak{e}_6=\Der( {\CD}(H_4(\Q)))$  by the group $  \mathbb{Z}_4^3 $, of type $( 48,15  )$;
  \item[\textbf{(7g11)}]    $\mathfrak{e}_7=\mathfrak{Instr}( {\CD}(H_4(\Q)))$ by the group $\mathbb{Z}_4^3   \times\mathbb{Z}_2 $, of type $( 102,14,1  )$.
\end{description}\smallskip

The restriction of the  $\mathbb{Z}_4^3$-grading on ${\CD}(H_4(\Q))$ to its structurable subalgebra  ${\CD}(H_4(\K))$
provides a $\mathbb{Z}_4^2\times\mathbb{Z}_2$-grading on ${\CD}(H_4(\K))$, which of course   can be used to get gradings on:
\begin{description}
\item[\textbf{(7g12)}]    $\mathfrak{e}_7=\mathfrak{Kan}( {\CD}(H_4(\K)))$  by the group $\mathbb{Z} \times \mathbb{Z}_4^2\times\mathbb{Z}_2 $, of type $(  67,27,4 )$;
  \item[\textbf{(7g13)}]    $\mathfrak{e}_7=\mathcal{U}( {\CD}(H_4(\K)))$ by the group $\mathbb{Z}_2^2 \times \mathbb{Z}_4^2\times\mathbb{Z}_2 =\mathbb{Z}_2^3 \times \mathbb{Z}_4^2$, of type $( 123,3,0,1  )$.
\end{description}\smallskip

\subsection{A fine $\mathbb{Z}_5^3$-grading}\label{subsec_elprimo5}

There is a $\mathbb{Z}_5^3$-grading on $\mathfrak{e}_8$ which seems not to be related with any of the previous constructions or structures. This grading appears in several contexts (for instance, \cite{Alek}  and \cite{Jordangradings}, and lately in \cite{dondelaZ52deE8}), due to its interesting properties: the zero homogeneous component is trivial (as in any fine grading by a finite group, see \cite[Corollary~5]{f4}) and all the other homogeneous components (in this case $124$) have the same dimension (so that in this case such dimension must be 2) and consist of semisimple elements (\cite[Lemma~1]{e6}). Moreover, given any $0\ne g\in\mathbb{Z}_5^3$, the subspace $\bigoplus_{i=1}^4\bigl(\mathfrak{e}_8\bigr)_{ig}$ is a Cartan subalgebra.
The following description can be found in \cite{Jordangradings}.

Let $V_1$ and $V_2$ be two vector spaces over $\FF$ of dimension 5 and let us consider the following $\mathbb{Z}_5$-graded vector space
$
L=L_{\bar 0}\oplus L_{\bar 1}\oplus L_{\bar 2}\oplus L_{\bar 3}\oplus L_{\bar 4},
$
for
\begin{equation}\label{eq_laZ5dee8graduacion}
\begin{array}{l}
L_{\bar 0}= \mathfrak{sl}(V_1)\oplus \mathfrak{sl}(V_2),  \\
L_{\bar 1}= V_1 \otimes \bigwedge^2V_2 , \\
L_{\bar 2}= \bigwedge^2V_1 \otimes \bigwedge^4V_2,  \\
L_{\bar 3}=  \bigwedge^3V_1 \otimes V_2 ,\\
L_{\bar 4}=   \bigwedge^4V_1 \otimes \bigwedge^3V_2.
\end{array}
\end{equation}
We can endow $L$ with a structure of $\mathbb{Z}_5$-graded Lie algebra, with the natural action of the semisimple algebra $L_{\bar 0}$ on each of the other homogeneous components. The brackets involving elements in different homogeneous components are given by suitable scalar multiples of the only $L_{\bar 0}$-invariant maps from $L_{\bar i}\times L_{\bar j}\to L_{\bar i+\bar j}$  (these scalars have been computed explicitly in \cite{poster}). The Lie algebra defined in this way is simple of dimension $248$, and hence  it provides a linear model of $\mathfrak{e}_8$. The philosophy of this kind of linear models can be found in \cite[Chapter~5, \S2]{enci}.

Let $\xi\in \FF$ be  a primitive fifth root of $1$, and take $B_1$ and $B_2$ bases of $V_1$ and $V_2$ respectively, and endomorphisms $b_1,c_1\in\End(V_1)$ and $b_2,c_2\in\End(V_2)$ whose coordinate matrices in the bases $B_i$ are
$$\begin{array}{ll}
b_1=P_\xi,&c_1=Q_5,\\
b_2=P_{\xi^2},&c_2=Q_5,
\end{array}
$$
defined as in Equation~(\ref{eq_Paulimatrices}).

Now the unique automorphisms $\Psi,\Psi'\in\Aut(L)$ whose restrictions to $L_{\bar1}$
are given by
$$
\begin{array}{l}
\Psi\vert_{L_{\bar1}}=b_1\otimes \wedge^2b_2,\\
\Psi'\vert_{L_{\bar1}}=c_1\otimes \wedge^2c_2,
\end{array}
$$
are order $5$ automorphisms which commute with the automorphism producing the $\mathbb{Z}_5$-grading on $L$ given by Equation~(\ref{eq_laZ5dee8graduacion}). Thus we obtained the desired grading by the group $\mathbb{Z}_5^3$, of type $(0,124)$.

\subsection{Gradings induced from other linear models}\label{subsec_explicacioneslinealesdeloscuatros}

We would like to explain a little bit the history of the search for the gradings described in Subsection~\ref{subsec_gradingsconcuatros}, which eventually lead to the quest for the grading  by $\mathbb{Z}_4^3$ on the simple structurable algebra of dimension $56$.

Consider the chain $E_6\subset E_7\subset E_8$ of exceptional groups.  The maximal abelian subgroup $\mathbb{Z}_4^3$ of $E_6$ is then also an abelian  subgroup of $E_7$ and also of $E_8$, predictably non-toral. This forced us to consider the order $4$ automorphisms of $\mathfrak{e}_7$ and  $\mathfrak{e}_8$.
First note that if we look at the subgroup $\mathbb{Z}_4^3$ of $E_6$, the three copies of $\mathbb{Z}_4$ involved  do not play the same role. One comes from $\theta$, an outer automorphism of $\mathfrak{e}_6$ producing the grading
\begin{equation}\label{eq_laZ4graddee6}
\mathfrak{e}_6=\left(\mathfrak{a}_3\oplus\mathfrak{sl}(V)\right)\oplus \left(V(2\lambda_1)\otimes V\right)\oplus \left(V(2\lambda_2)\otimes \FF\right)\oplus \left(V(2\lambda_3)\otimes V\right),
\end{equation}
for a two-dimensional vector space $V$, where the $\lambda_i$'s are the fundamental dominant weights for $\mathfrak{a}_3$. The other copies of $\mathbb{Z}_4$, restricted to  $\mathfrak{a}_3$, produce the $\mathbb{Z}_4^2$-grading obtained by means of Pauli matrices. More precisely, they correspond to the group
$$
\span{ (P_4,P_2),(Q_4,Q_2)}\subset \frac{\textrm{SL}(4)\times\textrm{SL}(2)}{\span{(\ii I_4,-I_2)}} \simeq \Cent_{E_6}\span{\theta}.
$$
The $\mathbb{Z}_4^3$-grading on $\mathfrak{e}_6$ is easily handled in this way, since we obtain concrete descriptions of the homogeneous components in terms of tensors of the  natural representations of $\mathfrak{sl}(4)$ and $\mathfrak{sl}(2)$.

Inspired by this, one can consider the automorphism of $\mathfrak{e}_7$ obtained by removing the black node of the extended Dynkin diagram (see \cite[Chapter~8]{libroKac})
\setlength{\unitlength}{0.04in}
\[
\begin{picture}(80,17)
	\put(10,12){\circle{2}}
	\put(20,12){\circle{2}}
    \put(30,12){\circle{2}}
    \put(40,12){\circle*{2}}
    \put(50,12){\circle{2}}
    \put(60,12){\circle{2}}
    \put(70,12){\circle{2}}
    \put(40,2){\circle{2}}
	
    \put(11,12){\line(1,0){8}}
    \put(21,12){\line(1,0){8}}
    \put(31,12){\line(1,0){8}}
    \put(41,12){\line(1,0){8}}
    \put(51,12){\line(1,0){8}}
    \put(61,12){\line(1,0){8}}
    \put(40,3){\line(0,1){8}}
	
	\put(9,14){\small 1}
	\put(19,14){\small 2}
    \put(29,14){\small 3}
    \put(39,14){\small 4}
    \put(49,14){\small 3}
    \put(59,14){\small 2}
    \put(69,14){\small 1}
    \put(42,1){\small 2}
	\end{picture}
\]
which produces a $\mathbb{Z}_4$-grading  $\mathcal{L}=\mathcal{L}_{\bar 0}\oplus\mathcal{L}_{\bar 1 }\oplus\mathcal{L}_{\bar2 }\oplus\mathcal{L}_{\bar 3}$ on the Lie algebra $\mathcal{L}\cong\mathfrak{e}_7$, where
$$
\begin{array}{l}
\mathcal{L}_{\bar 0}= \mathfrak{sl}(W_1)\oplus \mathfrak{sl}(V)\oplus \mathfrak{sl}(W_2),  \\
\mathcal{L}_{\bar 1}=  W_1\otimes V\otimes W_2 ,\\
\mathcal{L}_{\bar 2}=   \bigwedge^2W_1\otimes \FF\otimes \bigwedge^2W_2,\\
\mathcal{L}_{\bar 3}=   \bigwedge^3W_1\otimes V\otimes \bigwedge^3W_2,
\end{array}
$$
for $W_1=W_2$ and $V$ vector spaces of dimensions $4$ and $2$ respectively. (This gives the structure of the homogeneous components as modules for $\mathcal{L}_{\bar 0}$.)
We consider now the order $4$ automorphisms whose restrictions to $\mathcal{L}_{\bar1}$ are $P_4\otimes P_2\otimes P_4$,
$Q_4\otimes Q_2\otimes Q_4$, with $P_n$ and $Q_n$ defined as in     Equation~(\ref{eq_Paulimatrices}),
and the order $2$ automorphism determined by $w_1\otimes v\otimes w_2\mapsto w_2\otimes v\otimes w_1$. In this way, a $\mathbb{Z}_4^3\times \mathbb{Z}_2$-grading on $\mathfrak{e}_7$ is obtained (equivalent to \textbf{(7g11)}).

The same kind of arguments can be used to study $\mathbb{Z}_4$-gradings on $\mathfrak{e}_8$. Again,
  remove the black node of the extended Dynkin diagram of $\mathfrak{e}_8$
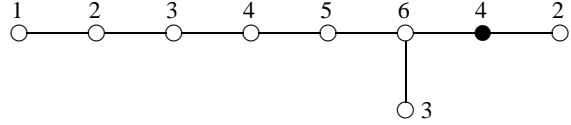
\setlength{\unitlength}{0.04in}
\[
\begin{picture}(80,20)
	\put(0,12){\circle{2}}
	\put(10,12){\circle{2}}
	\put(20,12){\circle{2}}
    \put(30,12){\circle{2}}
    \put(40,12){\circle{2}}
    \put(50,12){\circle{2}}
    \put(60,12){\circle*{2}}
    \put(70,12){\circle{2}}
    \put(50,2){\circle{2}}
	
    \put(1,12){\line(1,0){8}}
    \put(11,12){\line(1,0){8}}
    \put(21,12){\line(1,0){8}}
    \put(31,12){\line(1,0){8}}
    \put(41,12){\line(1,0){8}}
    \put(51,12){\line(1,0){8}}
    \put(61,12){\line(1,0){8}}
    \put(50,3){\line(0,1){8}}
	
	\put(-1,14){\small 1}
	\put(9,14){\small 2}
	\put(19,14){\small 3}
    \put(29,14){\small 4}
    \put(39,14){\small 5}
    \put(49,14){\small 6}
    \put(59,14){\small 4}
    \put(69,14){\small 2}
    \put(52,1){\small 3}
	
\end{picture}
\]
to get  an automorphism $\theta\in\Aut(\mathfrak{e}_8)$ producing a $\mathbb{Z}_4$-grading $\mathcal{L}=\mathcal{L}_{\bar 0}\oplus\mathcal{L}_{\bar 1 }\oplus\mathcal{L}_{\bar2 }\oplus\mathcal{L}_{\bar 3}$ on the Lie algebra $\mathcal{L}\cong\mathfrak{e}_8$ where,  as $\mathcal{L}_{\bar 0}$-modules, we have:
$$
\begin{array}{l}
\mathcal{L}_{\bar 0}= \mathfrak{sl}(W)\oplus \mathfrak{sl}(V),  \\
\mathcal{L}_{\bar 1}=  \bigwedge^2W\otimes V ,\\
\mathcal{L}_{\bar 2}=   \bigwedge^4W\otimes \FF,\\
\mathcal{L}_{\bar 3}=   \bigwedge^6W\otimes V,
\end{array}
$$
for $W$ and $V$ vector spaces of dimensions $8$ and $2$ respectively.
The centralizer can be checked to be
$$
\Cent_{E_8}(\theta) \simeq\frac{\PSL(8)\times \PSL(2)}{\span{(\ii I_8,-I_2)}},
$$
where the automorphism $\theta$ corresponds to the class $[(\xi I_8,I_2)]$ for $\xi\in \FF$ such that $\xi^2=\ii$. The MAD-groups of $E_8$ containing $\theta$ are MAD-groups of $\Cent_{E_8}(\theta)$. In particular, if we take $\theta$ together with
$$\begin{array}{ll}
{[}(I_2\otimes P_4,P_2)],\quad &[(P_2\otimes I_4,I_2)],\\
{[}(I_2\otimes Q_4,Q_2)],\quad &[(Q_2\otimes I_4,I_2)],
\end{array}
$$
we obtain   a $\mathbb{Z}_4^3\times \mathbb{Z}_2^2$-grading on $\mathfrak{e}_8$ (equivalent to \textbf{(8g12)});
and if we take it together with
$$
{[}(I_2\otimes P_4,P_2)],\quad {[}(I_2\otimes Q_4,Q_2)],\quad [(\left(\begin{smallmatrix}
\alpha &0 \\
0&\frac1\alpha \end{smallmatrix}\right)\otimes I_4,I_2)],
$$
we obtain the product of a one-dimensional torus and $\mathbb{Z}_4^3$, and hence a $\mathbb{Z}_4^3\times \mathbb{Z}$-grading on $\mathfrak{e}_8$ (equivalent to \textbf{(8g11)}).

According to \cite{paracasoinfinito},
this latter grading has to be induced by a fine grading with universal group $\mathbb{Z}_4^3$
on the   simple structurable algebra of dimension $56$. Unfortunately, the gradings on the structurable algebras are not yet classified, but in any case it was worth to find such a $\mathbb{Z}_4^3$-grading, because it lies behind several gradings on the simple Lie algebras $\mathfrak{e}_6$, $\mathfrak{e}_7$ and $\mathfrak{e}_8$ (those  in Subsection~\ref{subsec_gradingsconcuatros}). This was the starting point of \cite{Z43proceedings}.

However, the models based on the $\mathbb{Z}_4^3$-grading on the simple $56$-dimensional structurable algebra have a disadvantage over the linear models above, as the description of the homogeneous components involving pieces of $\mathfrak{trip}( {\CD}(H_4(C)))$ is not an easy task, and the type of the grading or the conjugacy classes of the automorphisms are neither easy to compute.

\subsection{Conclusion}

A large list of gradings on exceptional Lie algebras has been compiled here. All of them are fine and are described by their universal grading groups. We summarize them in the next result.

\begin{Theorem}\label{th:Main}
The following gradings on the simple Lie algebras of type $E$ are all fine:
\begin{itemize}
\item
The gradings on $\mathfrak{e}_6$ described as \textup{\textbf{(6gi)}}, $i=1,\dots,13$, whose universal groups are:
$\mathbb{Z}_2^6$, $\mathbb{Z}_3^4$, $\mathbb{Z}_2^3\times\mathbb{Z}_3^2$, $\mathbb{Z}_2\times\mathbb{Z}_3^3$,  $\mathbb{Z}^2\times\mathbb{Z}_3^2$, $\mathbb{Z} \times\mathbb{Z}_2^4$, $\mathbb{Z}^4 \times\mathbb{Z}_2$, $\mathbb{Z}^2 \times\mathbb{Z}_2^3$, $\mathbb{Z}  \times\mathbb{Z}_2^5$,  $\mathbb{Z}^2  \times\mathbb{Z}_2^3$,
 $ \mathbb{Z}_2^7$,  $\mathbb{Z}_4  \times\mathbb{Z}_2^4$,  $ \mathbb{Z}_4^3$.

\item  The gradings on $\mathfrak{e}_7$ described as \textup{\textbf{(7gi)}}, $i=1,\dots,13$,  whose universal groups are:
$\mathbb{Z}_2^7$,  $\mathbb{Z}_2^2\times\mathbb{Z}_3^3$, $\mathbb{Z} \times\mathbb{Z}_3^3$, $\mathbb{Z} \times\mathbb{Z}_2^5$, $\mathbb{Z}^4 \times\mathbb{Z}_2^2$, $\mathbb{Z}^3 \times\mathbb{Z}_2^3$,  $\mathbb{Z}  \times\mathbb{Z}_2^6$,  $\mathbb{Z}^2  \times\mathbb{Z}_2^4$,  $ \mathbb{Z}_2^8$,  $\mathbb{Z}_4  \times\mathbb{Z}_2^5$,  $\mathbb{Z}_4^3  \times\mathbb{Z}_2 $,  $\mathbb{Z}\times\mathbb{Z}_4^2  \times\mathbb{Z}_2 $, $\mathbb{Z}_2^3  \times\mathbb{Z}_4^2$.

\item  The gradings on $\mathfrak{e}_8$ described as \textup{\textbf{(8gi)}}, $i=1,\dots,13$,  whose  universal groups are:
$\mathbb{Z}_2^8$, $\mathbb{Z}_3^5$,  $\mathbb{Z}_6^3$, $\mathbb{Z}^2\times\mathbb{Z}_3^3$, $\mathbb{Z} \times\mathbb{Z}_2^6$, $\mathbb{Z}^4 \times\mathbb{Z}_2^3$,  $\mathbb{Z}  \times\mathbb{Z}_2^7$,  $\mathbb{Z}^2  \times\mathbb{Z}_2^5$,  $ \mathbb{Z}_2^9$,  $\mathbb{Z}_4  \times\mathbb{Z}_2^6$,  $\mathbb{Z}   \times\mathbb{Z}_4^3$,  $\mathbb{Z}_2^2  \times\mathbb{Z}_4^3$,    $\mathbb{Z}_5^3$.
    \end{itemize}
\end{Theorem}

\smallskip

\noindent \textbf{Conjecture:}  We think that these gradings exhaust the list of fine gradings, up to equivalence, on
$\mathfrak{e}_6$, $\mathfrak{e}_7$ and $\mathfrak{e}_8$, with the exception of the root space decompositions relative to a Cartan subalgebra.

If this conjecture were true, then there would be, up to equivalence,  $14$ fine gradings on each of the simple Lie algebras of type $E$. Therefore, there would be exactly $14$ conjugacy classes of maximal abelian diagonalizable subgroups of the algebraic group $\Aut(\mathfrak{e}_r)$, $r=6,7,8$.

This conjecture has been proved  for  $\mathfrak{e}_6$ in \cite{e6}. The cases of  $\mathfrak{e}_7$ and  $\mathfrak{e}_8$ remain open.
An strategy for the finite case is that every automorphism belonging to a MAD-subgroup of a connected and simply-connected group (like $E_8$) fixes a semisimple subalgebra, so that it corresponds to removing only one node in the extended Dynkin diagram. That fact implies that only a handful of automorphisms are possible and we are working in each case (see \cite{e8finito}) by studying the corresponding centralizers as in Subsection~\ref{subsec_explicacioneslinealesdeloscuatros}.
The group $\Aut(\mathfrak{e}_7)$ is not simply connected, 
but once one gets all the MAD-subgroups of $E_8$, much of the work is already done.
In order to deal with infinite MAD-groups, note that every grading by an infinite group is related to a grading by a root system, as proved in \cite{paracasoinfinito}. As these root-gradings are well known, the problem reduces to study some special gradings on the coordinate algebras. In many cases these are structurable algebras or related to them. Hence the problem is reduced to study gradings in algebras of relative low dimension (compared to the dimension of the exceptional simple Lie algebras).

\end{document}